\documentclass{article}

\usepackage{graphicx} 

\usepackage{amsthm,amssymb,amsmath,amscd,empheq,mathrsfs}
\usepackage{geometry}
\usepackage{todonotes}
\usepackage{comment}
\usepackage{amssymb}
\usepackage[shortlabels]{enumitem}
\usepackage{xcolor}
\usepackage[utf8]{inputenc}
\usepackage[T2A]{fontenc}
\usepackage{graphicx}
\usepackage{amssymb}
\usepackage{amsmath}
\usepackage{amssymb}
\usepackage{yhmath}

\usepackage{float}

\usepackage[style=alphabetic]{biblatex}

\usepackage{graphicx}
\usepackage{quiver}
\usepackage{tikz}
\usepackage{tikz-cd}
\usepackage{titlesec}\titleformat*{\section}{\LARGE\bfseries}
\usepackage{hyperref}
\usepackage{cleveref}

\hypersetup{
  colorlinks   = false,    
  urlcolor     = red,    
  linkcolor    = red,    
  citecolor    = red      
}

\usepackage[symbol]{footmisc}

\makeatletter
\newtheorem*{rep@theorem}{\rep@title}
\newcommand{\newreptheorem}[2]{%
\newenvironment{rep#1}[1]{%
 \def\rep@title{#2 \ref{##1}}%
 \begin{rep@theorem}}%
 {\end{rep@theorem}}}
\makeatother

\newreptheorem{theorem}{Theorem}

\theoremstyle{definition}
\newtheorem{example}{Example}[section]

\newtheorem{theorem}{Theorem}[section]
\newtheorem{lemma}[theorem]{Lemma}
\newtheorem{proposition}[theorem]{Proposition}
\newtheorem{corollary}{Corollary}[theorem]

\theoremstyle{definition}\newtheorem{definition}[theorem]{Definition}
\newtheorem{Problem}{Problem}[section]
\newtheorem{remark}[theorem]{Remark}

\numberwithin{equation}{theorem}

\DeclareMathOperator{\colim}{colim}

\DeclareMathOperator{\Spec}{Spec}
\DeclareMathOperator{\Spa}{Spa}
\DeclareMathOperator{\Shv}{Shv}
\DeclareMathOperator{\Spd}{Spd}

\DeclareMathOperator{\sld}{\square}
\DeclareMathOperator{\Sp}{Sp}

\DeclareMathOperator{\Fun}{Fun}

\DeclareMathOperator{\QCoh}{QCoh}

\DeclareMathOperator{\Mod}{Mod}
\DeclareMathOperator{\coMod}{coMod}

\DeclareMathOperator{\CAlg}{CAlg}
\DeclareMathOperator{\cCAlg}{cCAlg}
\DeclareMathOperator{\tri}{\triangleright}

\DeclareMathOperator{\SpecAn}{AnSpec}

\DeclareMathOperator{\AnR}{AnRing}
\DeclareMathOperator{\cyc}{cycl}
\DeclareMathOperator{\QC}{QCoh}
\DeclareMathOperator{\AnSt}{AnSt}

\DeclareMathOperator{\Cat}{Cat_{\infty}}
\DeclareMathOperator{\PrL}{Pr^{L}}

\DeclareMathOperator{\Z}{\mathbb{Z}}
\DeclareMathOperator{\Q}{\mathbb{Q}}

\DeclareMathOperator{\Div}{Div}
\DeclareMathOperator{\la}{la}
\DeclareMathOperator{\ct}{ct}

\DeclareMathOperator{\Sh}{Sh}

\DeclareMathOperator{\Ani}{Ani}

\DeclareMathOperator{\Aff}{Aff}

\DeclareMathOperator{\Hom}{{Hom}}

\title{Categorical K\"unneth formulas for analytic stacks}
\author{Youshua Kesting}

\begin{document}

\maketitle

\begin{abstract}
In \Cite{ben2010integral} the authors define a class of derived stacks, called "perfect stacks" and show that for this class the categories of quasi-coherent sheaves satisfy a categorical Künneth formula. Motivated to extend their results to the theory of analytic stacks as developed by Clausen-Scholze, we investigate categorical Künneth formulas for general $6$-functor formalisms.  As applications we show a general Tannakian reconstruction result for analytic stacks and, following recent work of Anschütz, Le Bras and Mann \Cite{anschütz20246functorformalismsolidquasicoherent}, show a $p$-adic version of Drinfeld's lemma for certain stacks that appear conjecturally in a categorical $p$-adic Langlands program.
\end{abstract}

\tableofcontents
\newpage
\section{Introduction}

The motivating problem of this article is to understand for which class of analytic stacks the categorical Künneth formula holds. Let us briefly recall the classical Künneth formula in algebraic geometry. Let $X\to \Spec(k),\ W \to \Spec(k)$ be two quasi-compact and separated schemes over a field $k$,denote by $D_{qcoh}(X)$ the derived category of $X$ with quasi-coherent cohomology sheaves and by $R\Gamma(X,-): D_{qcoh}(X) \to D_{qcoh}(k)$ the derived functor of global sections. Then the Künneth formula states that the canonical morphism 
\[
R\Gamma(X,\mathcal{O}_X)\otimes_k R\Gamma(W,\mathcal{O}_W)\to R\Gamma(X\times_{\Spec(k)} W,\mathcal{O}_{X\times_{\Spec(k)}W} )
\]
is an equivalence. In fact, this equivalence can even be observed on the level of the derived categories of $X$ and $W$ as has been investigated by Ben-Zvi, Francis and Nadler in the context of derived algebraic geometry. For $X\to Z, \ W \to Z$ any maps of perfect derived stacks and denoting by $\QC(X)$ the derived $\infty$-category of quasi coherent sheaves on $X$, they show (cf. \Cite[Theorem 4.14]{ben2010integral}) that one obtains a categorical Künneth formula
\begin{equation}\label{künnethforperfect}
    \QC(X)\otimes_{\QC(Z)}\QC(W) \cong \QC(X\times_Z W).
\end{equation}
Here the tensor product on the left hand side is the relative Lurie-tensor product in $\PrL$, the $\infty$-category of presentable $\infty$-categories. Let us recall that a derived stack $X$ is called perfect if its diagonal morphism $\Delta_X : X \to X\times X$ is affine, the $\infty$-category of quasi-coherent sheaves $\QC(X)$ is compactly generated and compact and dualisable objects coincide in $\QC(X)$ (cf. \Cite[Proposition 3.9]{ben2010integral}). Using the formalism of higher traces, they show furthermore that the existence of such categorical Künneth formulas implies many classical statements such as versions of the Atiyah-Bott-Lefschetz and the Riemann-Roch theorem (cf. \Cite[Theorem 1.4]{ben2013nonlinear}). \\
 The key idea of \Cite{ben2010integral} to prove \Cref{künnethforperfect} is to use that it is true for affine derived schemes, where it reduces to an abstract statement about $\infty$-categories of modules and to use descent of $\QC(-)$ for the general case. One difficulty arises due to the absence of a well-defined $6$-functor formalism for quasi-coherent sheaves, which they circumvent by showing proper base change and the projection formula for perfect stacks directly (cf. \Cite[Section 3.2]{ben2010integral}). \newline \\
 The later problem has been resolved by the recent theory of condensed mathematics and the theory of analytic stacks by Clausen-Scholze (cf. \Cite{clausen2022condensed},\Cite{AnSt}) which provides a vast framework to study both topology and algebraic, complex and $p$-adic geometry. In particular, they construct a $6$-functor formalism $D_{qc}(-)$ on the $\infty$-category $\AnSt$ of analytic stacks which satisfies a strong descent condition which they call $!$-descent.
Working in this generality, it becomes unnatural to impose the compact generatedness of $D_{qc}(-)$.  Moreover it become apparent in recent years that  there are many interesting $\infty$-categories such as the category of nuclear modules in condensed mathematics or categories of sheaves on locally compact Hausdorff spaces which are not compactly generated but which belong to a bigger class of presentable $\infty$-categories called dualisable categories.  The aim of this article is to show that the existence of a well-defined $6$-functor formalism with strong descent properties and the passage from compact-generated categories to dualisable categories yield strong generalizations of the results in \Cite{ben2010integral}. \\

More precisely, we can describe our problem as follows.

\begin{Problem}\label{motivationquestion}
For which morphisms $X\to Z,\ W\to Z$  of analytic stacks is the natural morphism
\begin{align}
 D_{qc}(X)\otimes_{D_{qc}(Z)} D_{qc}(W) \to D_{qc}(X\times_Z W)
\end{align}
an equivalence?    
\end{Problem}

Many of our arguments can be formulated for arbitrary $6$-functor formalisms. Motivated by finding a replacement of the notion of "perfect stacks" for general $3$- (or $6$)-functor formalisms, we adopt a maximalist approach and introduce the following definition.

\begin{definition}[\Cref{def.künneth6f}]
Let $D: Corr(C,E) \to \PrL$ be a $6$-functor formalism and $f: X\to Z$ be a map in $E$. We call $f$ Künneth if for all morphisms $W\to Z$ in $C$ the natural morphism 
\[D(X)\otimes_{D(Z)} D(W) \rightarrow D(X\times_Z W) \,, \, A\otimes B \mapsto p_X^*(A)\otimes p_W^*(B)
\]
is an equivalence, with $p_X : X \times_Z W \to X,\  p_W: X \times_Z W \to W$ the natural projection morphisms. We say $D$ satisfies Künneth for a class of morphisms $P\subset E$ if any  $f\in P$ is Künneth. We say $D$ satisfies Künneth if it satisfies Künneth for $P=E$.
\end{definition}
Note that for a general $6$-functor formalism there may be no non-trivial Künneth morphisms at all. For our applications to analytic stacks, we will be only concerned with the case where we start with a $3$-functor formalism $D_0: Corr(C_0,E_0) \to \PrL$ which already satisfies Künneth and consider its Kan extension to a $3$-functor formalism $D: Corr(C,E) \to \PrL$ on the category $C\coloneqq \Sh_{D_0}(C_0)$ of sheaves of anima on $C_0$ with respect to the $D_0$-topology, with $E$ an appropriate class of morphisms containing $E_0$ (see \Cite[Theorem 3.4.11]{heyer20246} for the precise conditions). \\
As we show, the class of Künneth morphisms enjoys many permanence properties (see \Cref{dualisbe stable under compposition}, \Cref{dualisable stable under bc}, \Cref{dualisability is !-local on source}, \Cref{dualisability is !-local on target}, \Cref{def open immersion for D Cop}, \Cref{open immersion and closed immersions are Künneth}):

\begin{proposition}\label{THM A}
Let $D: Corr(C,E) \to \PrL$ be a $6$-functor formalism.
\begin{enumerate}
\item (composition) If $f: Y \to Z$ and $g: X \to Y$ are Künneth, then $f\circ g: X \to Z$ is Künneth.
    \item (base change)  Let  $W \to Z$ be any morphism. If $f: Y \to Z$ is Künneth, then the base change $f': Y'\coloneqq Y\times_Z W \to W$ is Künneth.
   \item (open immersions) If $D: Corr(C,E) \to \Pr_{st}^L$ is stable and  $g: X \to Y$ in $E$ is a $D$-open immersion, then $f$ is Künneth.
    \item (locality on source) Let $g: X \to Y$ be any morphism in $E$. Let $(f_i: {X}_i\to X)_{i\in I}$ be a universal $!$-cover. If $g\circ f_{[n],i_\bullet} : X_{[n],i_\bullet}\coloneqq X_{i_0}\times_X ... \times_X X_{i_n} \to Y$ is Künneth for all $([n],i_\bullet)\in \Delta_I$, then $g$ is Künneth.
    \item  (locality on target) Let $f: X \to Z$ be any morphism in $E$ and $(g_i: Z_i\to Z)_{i\in I}$ be a universal $!$-cover. Denote by  $f_{([n],i_\bullet)}: X_{([n],i_\bullet)}\coloneqq X \times_Z Z_{[n],i_\bullet} \to Z_{[n],i_\bullet}$ the base change of $f$ by $g_{[n],i_\bullet}: Z_{[n],i_\bullet} \to Z$. If $f_{[n],i_\bullet}$ is Künneth for all $([n],i_\bullet)\in \Delta_I$, then $f$ is Künneth.
\end{enumerate}
\end{proposition}

Using these properties we show our main theorem (\Cref{Künnethextension-main corollary}, see \Cref{E' properties} for the notation). 

 \begin{theorem}\label{intro thmA} 
    Let $D: Corr(C,E) \to \PrL$ be a sheafy $6$-functor formalism on a subcanonical site $C$. Let $S\in C$ and assume that $D_S: ((C_E){{_{/S}}})^{op} \to \Mod_{D(S)}(\PrL)$ is symmetric monoidal. Then the extension $\tilde{D}_S: ((\Shv(C)_{\tilde{E}}){{_{/S}}})^{op}\to \Mod_{D(S)}(\PrL) $ is symmetric monoidal. The same is true if $S\in \Shv(C)$ has a representable diagonal $\Delta_S$ and a $!$-cover $U \to S$ with $U\in C$ and $D: Corr(C,E) \to \PrL$  is Künneth.
\end{theorem}

To formulate our first application, we introduce the following terminology (cf. \Cref{tannakian morphism}).
\begin{definition}
  Let $D: Corr(C,E)\to \CAlg(\PrL)$ be a $6$-functor formalism and $f: X \to S$ in $C$. We call $f$ Tannakian if the morphism induced by $D$
  \[
 \text{Hom}_{C_{/S}}(Y,X) \to \Fun_{D(S)}^{L,\otimes}(D(X),D(Y))
  \]
  is an equivalence for all $Y \in C_{/S}$.  We call $D$ Tannakian if all morphisms in $C$ are Tannakian.
\end{definition}

As a consequence of \Cref{intro thmA}, we prove the following Tannakian reconstruction result: 
\begin{theorem}\Cref{tannakian lift}
    Let $C$ be a subcanonical site and $D: Corr(C,E)\to \CAlg(\PrL)$ be a Tannakian sheafy $6$-functor formalism satisfying Künneth. Consider its extension \\ $D: Corr(\Shv(C),\tilde{E})\to \CAlg(\PrL)$. Let $f: X \to S$ be a morphism in $\tilde{E}$ and assume that there is a $!$-cover $g: S'\to S$ with $S'\in C$ and $g\in E_0$. Then $f$ is Tannakian.
\end{theorem}

We apply these results to the $6$-functor formalism  $\mathit{D}: Corr(\AnR, E) \to \PrL$ on the category $\AnR$ of analytic rings with $E$ the class of $!$-able maps (see  \Cref{$6$-functor formalism for analytic stacks} for a construction and the definition of $(D, \AnR,E)$). We show that $\mathit{D}: Corr(\AnR, E) \to \PrL$ is Künneth and Tannakian (\Cref{affinekünneth}) and obtain thus a partial answer to \Cref{motivationquestion}, using \Cref{intro thmA}.
\begin{corollary}\label{künethgeneral for analyticstacks}
Let $S\in \AnSt\coloneqq \Shv(\AnR))$ be an analytic stack with representable diagonal $\Delta_S$, and assume $S$ admits a $!$-cover $U \to S$ with $g: U \in \Aff$ an affine analytic stack and $g\in E_0$. Then any morphism of analytic stacks $X \to S$ in $\tilde{E}$ is Künneth and Tannakian.
\end{corollary}

Our interest in categorical Künneth formulas stems from finding an analogue of Drinfeld's lemma in the context of the $p$-adic local Langlands correspondence. Let us briefly recall the Drinfeld lemma in the $\ell$-adic setting (cf. \Cite[Chapter IV.7]{fargues2021geometrization}). Consider $E/\Q_p$ a finite extension, $k=\bar{\mathbb{F}}_q$ an algebraic closure of its residue field and $\Lambda=\bar{\Q}_{\ell}$. For $Y$ a small v-stack, consider the full subcategory $D_{dlb}(Y,\Lambda)\subset D_{et}(Y,\Lambda)$ of dualisable objects in the category of étale sheaves with coefficients in $\Lambda$ on $Y$. There is a small $v$-stack called $\Div_k^1$ whose significance comes from the fact that it "geometrises"  smooth $\Lambda$-representations of the Weil group $W_E$ of $E$  (cf. \Cite[Theorem V.1.1]{fargues2021geometrization}). Let $W$ be any small $v$-stack, the simplest instance of the Drinfeld lemma can then be phrased as an equivalence of categories
\[
D_{dlb}(\Div_k^1,\Lambda)\otimes_{D_{dlb}(\Spd(k),\Lambda)}D_{dlb}(W,\Lambda) \cong D_{dlb}(\Div_k^{1}\times_{\Spd(k)}W,\Lambda) 
\]
(cf. \Cite[Proposition 4.7.3]{fargues2021geometrization}). Note that the tensor product on the left hand side is the one induced from $\PrL$ on the category of idempotent complete stable $\infty$-categories.  \\
In the $p$-adic setting there is currently no conjecture along the lines of \Cite[Conjecture I.10.2]{fargues2021geometrization}.  Motivated by understanding the properties of pro-étale $\Q_p$-cohomology on rigid spaces, 
  Anschütz, Le Bras and Mann \Cite{anschütz20246functorformalismsolidquasicoherent} constructed recently a $6$-functor formalism $S \to D_{(0,\infty)}(S)\in \PrL$ on the category $\text{vStack}$ of small $v$-stacks, a variant of which might be a suitable category of coefficients for a $p$-adic Langlands correspondence. To support this claim, they introduce a small $v$-stack $\Div^1_E$ and show that it geometrises $(\varphi,\Gamma)$-modules on the Robba ring $\mathcal{R}_E$, where $\Gamma= \text{Gal}(E^{\cyc}/E)$. More precisely they show (\Cite[Proposition 6.3.15]{anschütz20246functorformalismsolidquasicoherent}) that the category of $(\varphi,\Gamma)$-modules on the $\mathcal{R}_E$ is contained fully faithfully in the category of dualisable objects in $D_{(0,\infty)}(\Div^1_E)$. A natural question from the context of a $p$-adic Langlands correspondence is thus whether the Drinfeld lemma holds in this setting, that is whether for $W$ any small $v$-stack we have an equivalence 
\[
 D_{(0,\infty)}({\Div}^1_{E}) \otimes_{D_{(0,\infty)}(\Spd(\mathbb{F}_q))} D_{(0,\infty)}(W) \cong D_{(0,\infty)}({\Div}^1_{E}\times_{\Spd(\mathbb{F}_q)}W).
\]
On $S=\Spa(R,R^+)$ an affinoid perfectoid space which admits a morphism of finite trg.dim to a totally disconnected perfectoid space, the category $D_{(0,\infty)}(S)\cong D_{\hat{\sld}}(Y_{(0,\infty,),S})$ is given by (cf. \Cite[Theorem 1.2.1]{anschütz20246functorformalismsolidquasicoherent}) the category of (modified) solid quasi-coherent sheaves on the analytic adic space
\[
Y_{(0,\infty),S}\coloneqq \Spa(W(R^+))\backslash V([\pi]p).
\]
Here $W(-)$ denotes the $p$-typical Witt vectors and $\pi\in R^+$ a pseudo-uniformiser. For simplicity, we will consider the case $E=\Q_p$ in the following but everything we say holds also for the a general finite extension $E$.

In order to avoid dealing with modified solid sheaves, we will consider the analytic stack $\Div_{\mathbb{Q}_p}^{1,\ct}\coloneqq Y_{(0,\infty),\Div^1_{\mathbb{Q}_p}}$ and the category $D(Y_{(0,\infty),\Div^1_{\mathbb{Q}_p}})$ instead of $D_{(0,\infty)}({\Div}^1_{\mathbb{Q}_p})$. By $v$-descent one can evaluate $D_{(0,\infty)}({\Div}^1_{\mathbb{Q}_p})$ by taking a  $v$-cover $S_{\bullet}\to {\Div}^1_{\mathbb{Q}_p}$, giving an equivalence  $D_{(0,\infty)}({\Div}^1_{\mathbb{Q}_p})\cong \lim_{[n]\in \Delta}D_{\hat{\sld}}(Y_{(0,\infty,),S_n})$. This identifies, up to the modification of the solid structure, with the descent datum  given by the cover $Y_{(0,\infty,),S_{\bullet}} \to \Div_{\mathbb{Q}_p}^{1,\ct}$ for $D_{qc}(-)$.  One technical reason we prefer working with analytic stacks directly is  the existence of $!$-descent for $D_{qc}(-)$. Another advantage of working directly on the level of analytic stacks is that one can  "geometrise" various kinds of $p$-adic representations within a single category. For example, for $G$ a $p$-adic Lie group, one has the analytic stacks $G$, $G^{\la}$  given by sending a compact open $U\subset G$ to the algebras $C(U,\Q_p)$, $C^{\la}(U,\Q_p)$ of continuous functions or locally analytic respectively. Consequently, the categories $D_{qc}(\SpecAn(\Q_p) /G^?)$, $? \in \{ \ct, \la\}$ correspond to certain continuous (respectively locally analytic) $\Q_p$-representations of $G$.   We will use this freedom to consider also a locally analytic version of the stack $\Div_{\mathbb{Q}_p}^{1,\ct}$ which we denote by $\Div_{\mathbb{Q}_p}^{1,\la}$ (cf. \Cref{Drinfeld cont und la}). By \Cref{künethgeneral for analyticstacks} the categorical Künneth formula holds for all of them. 
\begin{corollary}($p$-adic Drinfeld lemma)
    Let $? \in \{\ct,\la \}$. Then the morphism $\Div_{\Q_p}^{1,?}\to \SpecAn(\Q_p)$ is Künneth. In other words, for any analytic stack $W \to \SpecAn(\Q_p)$ we have equivalences
    \[
    \mathit{D}_{qc}(\Div_{\Q_p}^{1,?})\otimes_{\mathit{D}_{qc}(\Q_p)}  \mathit{D}_{qc}(W) \cong \mathit{D}_{qc}( \Div_{\Q_p}^{1,?} \times_{\SpecAn(\Q_p)}W).
    \]
\end{corollary}

The significance of Drinfeld's lemma in the $\ell$-adic setting of Fargues-Scholze lies in the construction of the spectral action (cf. \Cite[chapter VI, X]{fargues2021geometrization}). We hope that the above theorem will be useful to construct a spectral action in a $p$-adic (locally analytic) setting. \\

After this preprint was in its final form, Montagnani and Pavia informed us that they also independently obtained similar results to ours following similar ideas.

\subsection*{Outline}
In \Cref{Dualisable modules} we will collect some facts on presentable categories and the categories of modules over a presentably symmetric monoidal category. We will be brief and introduce just enough notation and definitions to discuss the characterisation of dualisable objects in $\Mod_E(\PrL)$ for $E\in \CAlg(\PrL)$ a commutative algebra in $\PrL$ due to Ramzi (cf.\Cite{ramzi2024dualizable}). We refer to \Cite[Chapter 1]{ramzi2024dualizable} for a more detailed discussion. \\
The main results of this paper are contained in \Cref{Künneth formulas for stacks}.  In \Cref{Künneth modules in $6$-functor formalisms} and \Cref{monoidal properties of Künneth} we introduce the notion of Künneth morphisms for arbitrary lax symmetric monoidal functors $D: C^{op} \to \PrL$ for a geometric setup $(C,E)$ and show first closure properties of Künneth maps.  In \Cref{Dualisable modules in $6$-functor formalisms} we recall the well known fact that for a $6$-functor formalism which satisfies Künneth, any morphism $X \to Y$ makes $D(X)\in \Mod_{D(Y)}(\PrL)$ into a dualisable module. In \Cref{Properties of Künneth morphisms} we recall the notion of $!$-descent and prove that the property of being Künneth is well behaved under descent. As an application, we prove a Tannakian lifting result in \Cref{Consequences: Tannakian lifting}. \\
Finally, in \Cref{$6$-functor formalism for analytic stacks}  we recall the $6$-functor formalism of analytic stacks and give our examples of Künneth morphisms coming from the Drinfeld lemma in \Cref{Drinfelds section}. Our main reference for $6$-functor formalisms is \Cite{heyer20246}.

\subsection*{Acknowledgments}
I want to thank Johannes Anschütz, Arthur-César Le Bras and Stefano Morra for our discussions concerning the Drinfeld lemma and for their guidance as advisors of my thesis, Cédric Pépin for his proposition to investigate the question of descent for Künneth morphisms and Adam Dauser for clarifying the role of representable morphisms and suggesting a counterexample for Künneth morphisms for analytic stacks. A special thanks goes to Greg Andreychev for reading an early draft of this paper and Rubis Bachelette for all the discussions around this paper.

\subsection{Notations and conventions}
\begin{itemize}
    \item We denote by $\PrL$ the (very large) $\infty$-category with objects presentable $\infty$-categories and morphisms colimit preserving functors. $\widehat{\Cat}$ denotes the $\infty$-category of large $\infty$-categories. We consider $\PrL$ as a symmetric monoidal $\infty$-category with respect to the Lurie-tensor product. 
    \item 
    We denote by $\Ani$ the $\infty$-category of spaces/anima and by $\Sp$ the $\infty$-category of spectra.

    \item For $C$ a symmetric monoidal $\infty$-category we denote by $\CAlg(C)$ the $\infty$-category of commutative algebra objects in $C$ and by $\cCAlg(C)\coloneq \CAlg(C^{op})^{op}$ the $\infty$-category of cocommutative coalgebra objects in $C$. For $A\in \CAlg(C), B\in \cCAlg(C)$ we denote by $\Mod_A(C)$ the $\infty$-category of $A$-modules and by $\coMod_B(C)$ the $\infty$-category of $B$-comodules in $C$. We denote by $\Pr_{st}^L\coloneq\Mod_{\Sp}(\PrL)$ the $\infty$-category of stable presentable $\infty$-categories.

    \item For $M,N \in \PrL$ we denote by $\Fun^L(M,N)$ the $\infty$-category of functors in $\PrL$. Similarly, for $E \in \CAlg(\PrL)$, $M,N \in \Mod_E(\PrL)$ we denote by  $\Fun_E^L(M,N)$ the $\infty$-category of functors in $\Mod_E(\PrL)$.  We denote the symmetric monoidal tensor product on $\Mod_E(\PrL)$ induced by the Lurie tensor product on $\PrL$ by $(-)\otimes_E (-)$ and for  $M,N \in \CAlg(\Mod_E(\PrL))$ we denote by $\Fun_E^{L,\otimes}(M,N)$ the mapping anima  in $\CAlg(\Mod_E(\PrL))$

    \item For $E \in \CAlg(\PrL)$ we denote by $\Cat_E$ the category of $E$-enriched $\infty$-categories. For $M\in \Cat_E$ an $E$-enriched category, we denote by $P_{E}(M)$ the category of $E$-enriched presheaves \Cite[Definition C.2.8]{heyer20246}.

    \item By "category" we will from now on mean $\infty$-category. We will refer to categories in which every inner horn has a unique filler as $1$-categories.

\end{itemize}

\newpage

\section{Dualisable modules}\label{Dualisable modules}
In this chapter we will recall the notion of presentable dualisable modules. Following \Cite[Chapter 1]{ramzi2024dualizable}, we introduce the notion of atomically generated categories in order to state an alternative characterisation of dualisability due to Ramzi \Cite[Proposition 1.40]{ramzi2024dualizable}. We will barely give any proofs and refer to \Cite[Chapter 1]{ramzi2024dualizable} for a more thourough treatment.
\begin{subsection}{Dualisable presentable categories}

Let us recall the notion of a dualisable object in a symmetric monoidal category. 

\begin{definition}
Let $C$ be a symmetric monoidal category. An object $P\in C$ is called dualisable if there are an object $P^\vee$, called the dual of $P$ and morphisms 
\[
ev : P\otimes P^\vee \to 1  \space \space , \space  \ co :1 \to P\otimes P^\vee
\]
called the evaluation and coevaluation maps such there are homotopy coherent diagrams 

\[\begin{tikzcd}
	P & {P\otimes P^\vee \otimes P} & {P^\vee} & {P^\vee \otimes P \otimes P^\vee} \\
	& P && {P^\vee}
	\arrow["{co \otimes id_P}", from=1-1, to=1-2]
	\arrow["{id_P}"', from=1-1, to=2-2]
	\arrow["{id_P \otimes ev}", from=1-2, to=2-2]
	\arrow["{id_{P^\vee}\otimes co}", from=1-3, to=1-4]
	\arrow["{id_{P^\vee}}"', from=1-3, to=2-4]
	\arrow["{ev \otimes id_{P^\vee}}", from=1-4, to=2-4]
\end{tikzcd}\]
We denote the full subcategory of $C$ spanned by the dualisable objects by  $C^{dlb}\subset C$.
\end{definition}
In the case of $C= \Pr_{st}^L$, there is the following description of dualisable objects due to Lurie.

\begin{theorem}\Cite[Proposition D.7.3.1]{lurie2018spectral}\label{dlb in PrL}
    Let $C\in \Pr_{st}^L$ be a stable presentable category. Then the following assertions are equivalent:
    \begin{enumerate}
\item $C$ is dualisable in $\Pr_{st}^L$.
\item $C$ is a retract of a compactly generated stable category.
\item There is an adjunction in $\Pr_{st}^L$
\begin{equation*}
\begin{tikzcd}
C \arrow[r,shift left=.5ex,"F"]
&
D  \arrow[l,shift left=.5ex,"G"]
\end{tikzcd}
\end{equation*}
with $D$ a compactly generated stable category such that $F$ is fully faithful and $G$ commutes with all colimits.
    \end{enumerate}
\end{theorem}

\subsection{Dualisable modules in presentable categories}
For the rest of this chapter let $E\in \CAlg(\PrL)$ be a commutative algebra object in $\PrL$. In the following subsection, we will be studying the category $\Mod_E(\PrL)$ and dualisable objects in it.  We will sometimes refer to an object of $\Mod_E(\PrL)$ as a presentable $E$-module.\\

As a first step to understand dualisable objects in $\Mod_E(\PrL)$, we want to obtain an analogue of the characterisation given in \Cref{dlb in PrL}, generalising from the case $E=\Sp$ to a general $E\in \CAlg(\PrL)$. To achieve that we need to generalise the notion of compactly generated categories to incorporate this additional module structure.
\begin{definition}\Cite[Definition 1.9]{ramzi2024dualizable}
    Let $f : M \to N$ be a functor in $\Mod_E(\PrL)$ with right adjoint $f_R$. We call $f$ an internal left adjoint if $f_R$ commutes with colimits and is $E$-linear, that is if the canonical projection map 
    \[
    x \otimes f_R(n) \to f_R(x \otimes n)
    \]
is an equivalence for all $n\in N$ and all $x\in E$.
\end{definition}
We denote by $\Mod(E)^{dlb}$ the (non full) subcategory of $\Mod_E(\PrL)$ spanned by dualisable objects and $1$-morphisms given by internal left adjoints. 
\begin{example}\Cite[Example 1.16]{ramzi2024dualizable}\label{examples of internal leftadjoint}
Let $C\in \PrL$ and $A \to B$ be a map in $\CAlg(C)$, then the functor given by base change of algebras  $\Mod_A(C)\to \Mod_B(C)$ is an internal left adjoint.
\end{example}

The following notion generalises the notion of compact objects in $\PrL$.
\begin{definition}\Cite[Definition 1.22]{ramzi2024dualizable}
 Let $M\in \Mod_E(\PrL)$. An object $x\in M$ is called $E$-atomic (or just atomic, if $E$ is clear from the context) if the functor $x \otimes(-) :E \to M$ is an internal left adjoint.  
\end{definition}
We see that if $E=\Sp$ is the category of spectra the $\Sp$-atomic objects in $\Mod_{\Sp}(\PrL)\cong \Pr_{st}^L$ are exactly the compact objects. Following \Cite{ramzi2024dualizable}, we now introduce the following generalization of compactly generated categories.

\begin{definition}\Cite[Definition 1.27]{ramzi2024dualizable}
    An $E$-module $M\in \Mod_E(\PrL)$ is called $E$-atomically generated if the smallest full sub-$E$-module of $M$ closed under colimits and containing the atomics of $M$ is $M$ itself.
\end{definition}
In the following we will make use of the notion of $E$-enriched categories, our main reference is \Cite{heine2023equivalence} and \Cite[Appendix C]{heyer20246}. For $M\in \Cat_E$ an $E$-enriched category, we denote by  $\widehat{\Cat}_E$  the category of large $E$-enriched categories. By \Cite[Theorem 1.2]{heine2023equivalence} there is a natural functor 
\[
\mu : \Mod_E(\PrL) \to  \widehat{\Cat}_E
\]
which is fully faithful on hom categories. 
\begin{lemma}\Cite[Observation 1.28]{ramzi2024dualizable}
An $E$-module $M\in \Mod_E(\PrL)$ is $E$-atomically generated  if and only if it is equivalent to a category of $E$-enriched presheafs  $P_{E}(M_0)$ for $M_0$ some $E$-enriched category.
\end{lemma}
$E$-atomically generated categories are thus generalisations of compactly generated categories. One crucial step in the proof of \Cref{dlb in PrL} is to show that compactly generated categories are dualisable in $\Pr_{st}^L$. The following Lemma due to Berman \Cite[Theorem 1.7]{berman2020enriched} shows that the analogous assertion is true for $E$-atomically generated modules in $\Mod_E(\PrL)$.
\begin{lemma}\Cite[Proposition 1.40]{ramzi2024dualizable}
 Any $E$-atomically generated module $M\in \Mod_E(\PrL)$ is dualisable in $\Mod_E(\PrL)$.
\end{lemma}

We now have all the necessary concepts to obtain a characterisation of dualisable objects in $\Mod_E(\PrL)$.
\begin{theorem}\Cite[Theorem 1.49]{ramzi2024dualizable}\label{dlb relative in PrL}
    Let $E\in \CAlg(\PrL)$ and  $C\in \Mod_E(\PrL)$ be a presentable $E$-module. Then the following assertions are equivalent:
    \begin{enumerate}
\item $C$ is dualisable in $\Mod_E(\PrL)$.
\item $C$ is a retract of an  $E$-atomically generated category.
\item There is an internally left adjoint fully faithful embedding $j : C \hookrightarrow D$ for some $E$-atomically generated category $D$.
    \end{enumerate}
\end{theorem}

\subsection{Extension and restriction of modules}
We end this chapter with a few lemmas on the behaviour of dualisable modules under extensions and restrictions of modules. 
\begin{lemma}\label{tensor in Prl commutes with colimits}
 Any morphism $g: V \to E$ in $\CAlg(\PrL)$ induces a functor $G_L\coloneq E \otimes_V (-): \Mod_V(\PrL) \to  \Mod_E(\PrL)$ which is symmetric monoidal and commutes with small colimits.  
\end{lemma}
\begin{proof}
    The existence and symmetric monoidality of the functor $G_L$ as a functor in $\widehat{\Cat}$ is \Cite[Remark 4.5.3.2]{lurie2017higher}. We note that the cocartesian fibration $\Mod(\PrL) \to \CAlg(\PrL) \times  N(Fin_\ast)$ is presentable, in the sense that its classifying functor $\CAlg(\PrL) \times  N(Fin_\ast) \to \widehat{\Cat} $ factorises over $\PrL \subset \widehat{\Cat}$. The fact that $G_L$ commutes with small colimits follows since $ \otimes: \PrL \times \PrL \to \PrL $ commutes with small colimits in each variable separately.
\end{proof}

\begin{lemma}\label{tensor with dualisable object commmutes lim}
Let $g: V \to E$ be a morphism in $\CAlg(\PrL)$. Assume $E$ is dualisable in $\Mod_V(\PrL)$. Then $G_L : \Mod_V(\PrL) \to  \Mod_E(\PrL)$ commutes with limits.
\end{lemma}
\begin{proof}
Let $F: I \to \Mod_V(\PrL)$ be a limit diagram, then we have equivalences 
\begin{align*}
  E \otimes_V \lim_n F(n) &\cong \Fun^L_V(V, E \otimes_V \lim_n F(n)) \\
  &\cong \Fun^L_V(E^\vee, \lim_n F(n)) \\
  &\cong \lim_n\Fun^L_V(E^\vee, F(n)) \\
   & \cong \lim_n E \otimes_V F(n).  
\end{align*}
\end{proof}

\begin{lemma}\label{extensionrestriction}
Let $g: V \to E$ be a morphism in $\CAlg(\PrL)$. The right adjoint $G_R$ of $G_L$ preserves dualisable objects if and only if $E$ is dualisable in $\Mod_V(\PrL)$. It reflects dualisable objects if $E$ is dualisable as a $E\otimes_V E$-module.
\end{lemma}
\begin{proof}
The claim that $G_R$ preserves dualisable objects if and only if $E$ is dualisable in $\Mod_V(\PrL)$ follows by \Cite[Proposition 4.6.4.4. (8)]{lurie2017higher}. The second part follows by \Cite[Proposition 4.6.4.12 (6)]{lurie2017higher}.
\end{proof}

\end{subsection}

\newpage

\section{Künneth formulas for stacks}\label{Künneth formulas for stacks}
We introduce the notion of Künneth morphisms for general $3$-functor formalisms and discuss its relation to the dualisable presentable modules discussed in the previous chapter. For a $6$-functor formalism $D: Corr(C,E)\to \PrL$, a morphism $f: X \to Z$ in $E$ will be called Künneth if a general categorical Künneth formula holds, by which we mean that for all $W\to Z$ in $C$ the natural morphism
\[
m_{X,W}: D(X) \otimes_{D(Z)} D(W) \to D(X\times_Z W)
\]
is an equivalence. We show that Künneth morphisms are stable under composition and base change (\Cref{dualisbe stable under compposition}, \Cref{dualisable stable under bc}) and in certain cases,  $!$-local on the source and target (\Cref{dualisability is !-local on source}, \Cref{dualisability is !-local on target}). Finally, we recall the $6$-functor formalisms of analytic stacks and give examples of Künneth morphisms in this context, giving a generalisation of the discussion in \Cite[Chapter 3-4]{ben2010integral} in the context of quasi-coherent sheaves on perfect stacks. We refer to \Cite[Chapter 1-4]{heyer20246} for a reference on general $6$-functor formalisms.
\subsection{Künneth morphisms}\label{Künneth modules in $6$-functor formalisms}
A very powerful tool to construct $3$-functor formalisms is the construction result of Heyer-Mann \Cite[Proposition 3.3.3]{heyer20246}. We start by recalling the following definition.
\begin{definition}\Cite[Definition 2.1.1, Remark 2.1.2]{heyer20246}
A geometric setup is a tuple $(C,E)$ of a category $C$ together with a homotopy class of edges $E$ in $C$ such that the following conditions are satisfied.
\begin{enumerate}
    \item $E$ contains all isomorphisms. 
    \item $E$ is stable under compositions and under pullback along edges in $C$. 
    \item If $f :X \to Y$ lies in $E$, then so does its diagonal $\Delta_f: X \to X \times_Y X$.
\end{enumerate}
  
\end{definition}
For the construction result of Heyer-Mann, assume we are given a geometric setup $(C,E)$, where $C$ is a category with finite limits, and a lax symmetric monoidal functor 
\begin{align*}
  D : C^{op} \to \CAlg(\PrL), \ X &\mapsto D(X)   \\
   (f: Y \to X) &\mapsto  D(f)\coloneqq f^* : D(X) \to D(Y)
\end{align*}
where we consider $C$ with its cartesian monoidal structure. We will furthermore assume that the class $E$ decomposes into a suitable class $P\subset E$ of "proper maps" and a class $I\subset E$ of "local isomorphisms, such that any $f\in E$ can be written as a composition of a "local isomorphism" and a "proper map". We now want to extend the functor $D$ to a $3$-functor formalism 
\[
D: Corr(C,E) \to \CAlg(\PrL)
\]
such that for any $f\in P$ the functor $f_!$ is a right adjoint of $f^*$ and for any $f\in I$ the functor $f_!$ is a left adjoint to $f^*$. The construction theorem of Mann states that such an extension to a $3$-functor formalism exists given a certain amount of data (cf. \Cite[Proposition 3.3.3]{heyer20246} for the conditions).\\

In view of this construction theorem, we will begin by discussing categorical Künneth formulas quite abstractly. Let $C$ be a category with finite limits and consider  $C$ with its cartesian monoidal structure. Let  $D: C^{op}\to \CAlg(\PrL)$ be a lax symmetric monoidal functor. Let us recall how to define the morphisms $m_{X,W}$ in this context.
Let $X,W\in C$ be two objects and $p_X : X\times_Z W \to X$ and $p_W : X \times_Z W\to W$ the canonical projections. Then the functors  $p_X^*$ and $p_W^*$ are $D(Z)$-linear and combine to a $D(Z)$-bilinear functor 
\[
D(X)\times D(W) \to D(X \times_Z W), \quad (A, B) \mapsto p_X^*(A) \otimes p_W^*(B).
\]
By the universal property of the Lurie-tensor product, this functor extends to a $D(Z)$-linear functor which we denote by
\[
m_{X,W}: D(X)\otimes_{D(Z)} D(W) \to D(X \times_Z W).
\]

\begin{lemma}\label{naturality of m1}
Let $(C,E)$ be a geometric setup, and consider $C^{op}$ as a symmetric monoidal category with its coCartesian symmetric monoidal structure and $D: C^{op} \to \CAlg(\PrL)$ be a lax symmetric monoidal functor. Let $X\to Z, W \to Z$ and $f: Y\to W$ be any maps in $C$. Then the following diagram commutes
\[\begin{tikzcd}
	{D(X)\otimes_{D(Z)}D(W)} & {D(X)\otimes_{D(Z)}D(Y)} \\
	{D(X\times_Z W)} & {D(X\times_Z Y)}
	\arrow["{id_{D(X)}\otimes f^*}", from=1-1, to=1-2]
	\arrow["{(id_X\times f)^*}", from=2-1, to=2-2]
	\arrow["{m_{X,W}}"', from=1-1, to=2-1]
	\arrow["{m_{X,Y}}", from=1-2, to=2-2]
\end{tikzcd}
\]
\end{lemma}
\begin{proof}
This follows from the lax monoidality of $D: C^{op} \to \CAlg(\PrL)$.
\end{proof}

We introduce the following definition.
\begin{definition}\label{def.künneth6f}
Let $(C,E)$ be a geometric setup, and consider $C^{op}$ as a symmetric monoidal category with its coCartesian symmetric monoidal structure. 
\begin{enumerate}
    \item Let $D: C^{op} \to \CAlg(\PrL)$ be a lax symmetric monoidal functor. Let $f: X\to Z$ be a map in $E$. We call $f$ Künneth if for all objects $W \in C_{/Z}$ the natural morphism 
\[
m_{X,Y} : D(X)\otimes_{D(Z)} D(W) \to  D(X\times_Z W)
\]
is an isomorphism. We say that $D$ satisfies Künneth for a class $P \subset E$ if any morphism $f \in P$ is Künneth. We say $D$ satisfies Künneth if it satisfies Künneth for $P=E$.
\item If $D: Corr(C,E) \to \CAlg(\PrL)$ is a $3$-functor formalism and $f: X \to Y$ in $E$, we call $f$ Künneth if it is Künneth for the functor $C^{op}\to Corr(C,E)\to \CAlg(\PrL)$. 
\end{enumerate}
\end{definition}

\subsection{Monoidal properties of Künneth morphisms}\label{monoidal properties of Künneth}
In this subsection we consider Künneth morphisms for $D: C^{op} \to \CAlg(\PrL)$  a lax symmetric monoidal functor for $(C,E)$ a geometric setup. We show certain closure properties of Künneth morphisms such as stability under composition and base change and give first examples of Künneth morphisms. 

\begin{lemma}\label{dualisbe stable under compposition}
Let $D: C^{op} \to \CAlg(\PrL)$ be a lax symmetric monoidal functor, $f: X \to Y$, $g: Y \to Z$ two morphisms. If $f$ and $g$ are Künneth, then $g\circ f$ is Künneth. 
\end{lemma}
\begin{proof}
 Let $W \to Z$ be any morphism in $C$. Since $f$ is Künneth and using the isomorphism $X \times_Z W \cong X \times_Y (Y\times_Z W)$ we obtain an equivalence 
\[
m_{X,Y\times_ZW}: D(X) \otimes_{D(Y)}D(Y\times_Z W)\cong D(X\times_ZW ).
\]
Using that $g$ is Künneth, we obtain equivalences 
\[
D(X) \otimes_{D(Y)}D(Y\times_Z W)\cong D(X)\otimes_{D(Z)}D(W).
\]
By functoriality of pullbacks and \Cref{naturality of m1} we see that the resulting equivalence $D(X)\otimes_{D(Z)}D(W) \to D(X\times_Z W)$ is given by $m_{X,W}$ which proves that $g\circ f$ is Künneth. 
\end{proof}
\begin{lemma}\label{dualisable stable under bc}
Let $D: C^{op} \to \CAlg(\PrL)$ be a lax symmetric monoidal functor, $f: X \to Z$, $g: Y \to Z$ two morphisms. If $f$ is Künneth, then the base change $f' : X' \coloneq X\times_Z Y \to Y$ of $f$ along $g$ is Künneth.
\end{lemma}
\begin{proof}
 Let $W \to Y$ be an arbitrary morphism. We have a chain of equivalences 
 \begin{align*}
     D(X') \otimes_{D(Y)}D(W) & \cong (D(X)\otimes_{D(Z)}D(Y) )\otimes_{D(Y)}D(W) \\
     & \cong D(X) \otimes_{D(Z)}D(W) \\
     & \cong D(X' \times_Y W)
 \end{align*}
 Here we used that $f$ is Künneth in the first line, the natural isomorphism in the second line and the natural isomorphism $X' \times_Y W \cong X\times_Z W$ together with the fact that $f$ is Künneth in the last line.
\end{proof}

\begin{lemma}\label{dualisbe and cancelation}
Let $D: C^{op} \to \CAlg(\PrL)$ be a lax symmetric monoidal functor, $f: X \to Y$, $g: Y \to Z$ two morphisms. Assume that  $g$ and the diagonal $\Delta_g: Y \to Y\times_Z Y$ are Künneth. If  $g\circ f$ is Künneth, then $f$ is Künneth.
\end{lemma}
\begin{proof}
We write $f$ as the composition $ X \overset{(id,f)}{\to} X \times_Z Y \overset{pr_2}{\to} Y$. By \Cref{dualisbe stable under compposition}, it suffices to show that $(id,f)$ and $pr_2$ are Künneth. The morphism $(id,f)$ is given by the base change of $(f,id): X\times_Z Y \to Y \times_Z Y$  along the diagonal $\Delta_g : Y \to Y\times_Z Y$ an is thus Künneth by \Cref{dualisable stable under bc}. Similarly, the morphism $pr_2$ is the base change of $g \circ f$ along $g$ and is thus Künneth by \Cref{dualisable stable under bc}.
\end{proof}

\begin{lemma}\label{dualisability is stable under product}
Let  $D: C^{op} \to \CAlg(\PrL)$ be a lax symmetric monoidal functor. If  $f: X_1 \to Z$, $g: X_2 \to Z$ are Künneth, then $X_1\times_Z X_2 \to Z$ is Künneth.
\end{lemma}
\begin{proof}
This follows from \Cref{dualisable stable under bc} and \Cref{dualisbe stable under compposition}.
\end{proof}

We will now give first examples of Künneth morphisms. We recall the following definition.
\begin{definition}\Cite[Lemma 6.4, Proposition 6.5]{clausen2022condensed}\label{abstractdef of open and closed in Sym}
Let $F : \mathcal{D} \to \mathcal{C}$ be a map in $\CAlg(\Pr_{st}^L)$.  \\
We call $F$ a categorical open immersion, if $F$ has a fully faithful left adjoint $F_L: \mathcal{C} \to \mathcal{D}$ which satisfies the projection formula, i.e.  for any $M\in \mathcal{C}$ and any $N\in \mathcal{D}$ the natural map 
    \[
   F_L(M\otimes F(N)) \to  F_L(M) \otimes N
    \]
    is an equivalence. 

\end{definition}
\begin{proposition}\Cite[Proposition 6.5]{clausen2022condensed}\label{alternative defclosed open inSym}
Let $F : \mathcal{D} \to \mathcal{C}$ be a map in $\CAlg(\Pr_{st}^L)$. Then $F$ is a categorical open immersion if and only if there is a (necessarily unique) idempotent algebra $A\in \mathcal{D}$ such that $F(A)\cong 0$ and the induced map $\mathcal{D}/\Mod_A\mathcal{D} \to \mathcal{C}$ is an equivalence.

\end{proposition}
\begin{proof}
  Assume that $F$ is a categorical open immersion. We claim that $A\coloneqq \mathrm{cofib}(F_L(1) \to 1)\in \mathcal{D}$ is an idempotent algebra. Indeed, note that $F(A)\cong \mathrm{cofib}(F (F_L(1))\to 1)\cong 0$ since $F_L$ is fully faithful.  We thus have equivalences
\begin{align*}\label{idempotetn algebra etc}
A\otimes A & \cong \mathrm{cofib}(F_L(1)\otimes A \to A)\\
& \cong \mathrm{cofib}(F_L(F(A)) \to A) \\
& \cong \mathrm{cofib}(0 \to A)\cong A
\end{align*}
where the first equivalence follows since cofibers commute with tensor products and the second line follows by the projection formula and the last line follows since $F(A)\cong 0$. Thus $A\in \mathcal{D}$ is an idempotent algebra and it suffices to show that
\[
\mathrm{ker}(F)\cong \Mod_A\mathcal{D}
\]
as full subcategories of $\mathcal{D}$. If $M\in \mathrm{ker}(F)$, then $F_L(1)\otimes M \cong F_L (F(M)) \cong 0$ by the projection formula, so $A \otimes M\cong M $ and thus $M$ lies in $\Mod_A\mathcal{D}$. If $M\in \Mod_A\mathcal{D}$ we have $A \otimes M \cong M$ and thus $F(A \otimes M) \cong F(A) \otimes F(M) \cong 0$ since $F(A)=0$. Conversely, assume that there exists an idempotent algebra $A\in \mathcal{D}$ such that $F :\mathcal{D} \to \mathcal{D}/\Mod_A\mathcal{D} \cong \mathcal{C}$. Then $F(-)\cong \mathrm{fib}(1\to A)\otimes (-) $ and  $F_L$ is the inclusion, so the projection formula is trivially satisfied.
\end{proof}
\begin{lemma}\Cite[Lemma 6.4, Corollary 6.6]{clausen2022condensed}\label{compos bc and cancel for openclosedimmsersion in sym}
\begin{enumerate}
    \item Let $F : \mathcal{D} \to \mathcal{C}$, $G:\mathcal{C}\to \mathcal{B} $ be categorical open immersions, then $G\circ F$ is a categorical open immersion.
\item Let $F : \mathcal{D} \to \mathcal{C}$ be a categorical open immersion, $G: \mathcal{D} \to \mathcal{B}$ be a functor in $\CAlg(\Pr_{st}^L)$. Then the pushout $\tilde{F}: \mathcal{B} \to \mathcal{C}\otimes_{\mathcal{D}}\mathcal{B}$ is again a categorical open immersion. 
\item Let $F : \mathcal{D} \to \mathcal{C}$, $G:\mathcal{C}\to \mathcal{B} $ be maps in $\CAlg(\Pr_{st}^L)$. If $G\circ F$ is a categorical open immersion, then $G$ is a categorical open immersion.
    \item Let $G : \mathcal{D} \to \mathcal{B}$ be a map in $\CAlg(\Pr_{st}^L)$ and consider the pushout diagram in $\CAlg(\Pr_{st}^L)$.
    \[\begin{tikzcd}
	{\mathcal{B}} & {\mathcal{B}\otimes_{\mathcal{D}}\mathcal{C}} \\
	{\mathcal{D}} & {\mathcal{C}}
	\arrow["{\tilde{F}}", from=1-1, to=1-2]
	\arrow["{G}"', from=2-1, to=1-1]
	\arrow["F", from=2-1, to=2-2]
	\arrow["\tilde{G}"', from=2-2, to=1-2]
\end{tikzcd}\]
If  $F$ is a categorical open immersion, the diagram above is left adjointable.
\end{enumerate}
\end{lemma}

\begin{definition}\label{def open immersion for D Cop}
$D: C^{op} \to \CAlg(\Pr_{st}^L)$ be a lax symmetric monoidal functor and $f: X \to Y$ any morphism in $E$. We call $f$ a $D$-closed  (respectively $D$-open) immersion if $D(f)\coloneqq f^* : D(Y) \to D(X)$ is a categorical closed (respectively open) immersion in the sense of \Cref{abstractdef of open and closed in Sym}.
\end{definition}
We will now see that in the context of a $3$-functor formalisms, these notions yield our first examples of Künneth morphisms.
\begin{proposition}\label{open immersion and closed immersions are Künneth}
$D: Corr(C,E) \to \CAlg(\Pr_{st}^L)$ be a $3$-functor formalism and $f: X \to Y$ any morphism in $E$. 
Let $f$ be a $D$-open immersion, then $f$ is Künneth.

\end{proposition}
\begin{proof}
 Let $f: X \to Y$ be a $D$-open immersion and $g: W \to Y$ any morphism in $C$. Let $I \in D(Y)$ be the idempotent algebra associated with $f$ via \Cref{alternative char. of open and proper in aff}. We have a localisation sequence
 \[
 \Mod_ID(Y) \hookrightarrow D(Y) \overset{f^*}{\to} D(X).
 \]
 Now the Lurie-tensor product in $\Mod_{D(Y)}(\PrL)$ preserves localizations (cf.\Cite[Corollary 1.46]{ramzi2024dualizable}), thus after applying $(-)\otimes_{D(Y)}D(W)$ we obtain a localization sequence
\[\begin{tikzcd}
	{\Mod_ID(Y)\otimes_{D(Y)}D(W)} & {D(Y)\otimes_{D(Y)}D(W)} & {D(X)\otimes_{D(Y)}D(W)} \\
	{\Mod_{g^*(I)}D(W)} & {D(W)} & {D(X\times_YW)} \\
	{}
	\arrow[hook, from=1-1, to=1-2]
	\arrow["\wr"', from=1-1, to=2-1]
	\arrow["{f^*\otimes id_{D(W)}}", shift left, from=1-2, to=1-3]
	\arrow["{\wr \ m_{Y,W}}", from=1-2, to=2-2]
	\arrow["{m_{X,W}}", from=1-3, to=2-3]
	\arrow[hook, from=2-1, to=2-2]
	\arrow["{\tilde{f}^*}", from=2-2, to=2-3]
\end{tikzcd}\]
Note that $g^*(I)$ is idempotent so the lower left functor is again fully faithful. The right square commutes by \Cref{naturality of m1}. We need to show $\text{ker}(\tilde{f}^*)\cong \Mod_{g^*(I)}D(W)$ as full subcategories of $D(W)$. The inclusion $\Mod_{g^*(I)}D(W)\subset \text{ker}(\tilde{f}^*)$ is obvious by the commutativity of the diagram. Now let $M \in  \text{ker}(\tilde{f}^*)$, we want to show that $g^*(I) \otimes M \cong M$. This follows by the following calculation
\begin{align*}
    g^*(I) \otimes M & \cong \text{cofib}(g^*f_!(1)\otimes M \to M) \\
    & \cong \text{cofib}(\tilde{f}_!(1)\otimes M \to M) \\
    & \cong \text{cofib}(\tilde{f}_!(\tilde{f}^*(M)) \to M) \\
    & \cong  \text{cofib}(0\to M)\cong M.
\end{align*}
Thus we conclude that $\tilde{f}$ is a $D$-open immersion and since $m_{Y,W}$ and the left vertical morphism are isomorphisms also $m_{X,W}$ is an isomorphism.
\end{proof}

\begin{corollary}\label{bc for open immersions}
 $D: Corr(C,E) \to \CAlg(\PrL)$ be a $3$-functor formalism and  $f: X \to Y$ a $D$-open immersion, $g: Z \to Y$ any morphism in $C$. Then the base change $\tilde{f}: X \times_Y Z \to Z$ is again a $D$-open immersion.
\end{corollary}
\begin{proof}
    This follows by \Cref{open immersion and closed immersions are Künneth} together with point $2$ in \Cref{compos bc and cancel for openclosedimmsersion in sym}.
\end{proof}

\subsection{Dualisable modules in $6$-functor formalisms}\label{Dualisable modules in $6$-functor formalisms}
Now let $D : Corr(C,E) \to \PrL$ be a $6$-functor formalism satisfying Künneth and $f: X \to Z$ a map in $E$. We want to see that $D(X)$ is then automatically dualisable in $\Mod_{D(Z)}\PrL$ via the map $f^*$.

\begin{lemma}\label{funnydiagrams}
Let $(C,E)$ be a geometric setup,  $\pi : X \to Z$ a map in $E$ and $\Delta : X \to X\times_Z X$ the diagonal map. Then the following two compositions  in $\mathrm{Corr}(C,E)$ represent the identity. In particular, $X$ is self-dual in $\mathrm{Corr}(C,E)$.
\[\begin{tikzcd}
	& {X\times_Z X} & {} & {X\times_Z X} \\
	X && {X\times_Z X \times_Z X} & {} & X
	\arrow["{\Delta\times id_X}", from=1-2, to=2-3]
	\arrow["{\pi\times id_X}"', from=1-2, to=2-1]
	\arrow["{id_X \times \Delta}"', from=1-4, to=2-3]
	\arrow["{id_X \times \pi}", from=1-4, to=2-5]
\end{tikzcd}\]
\[\begin{tikzcd}
	& {X\times_Z X} & {} & {X\times_Z X} \\
	X && {X \times_Z X \times_Z X} & {} & X.
	\arrow["{id_X \times \Delta}", from=1-2, to=2-3]
	\arrow["{\Delta \times id_X }"', from=1-4, to=2-3]
	\arrow["{\pi \times id_X}", from=1-4, to=2-5]
	\arrow["{id_X \times \pi}"', from=1-2, to=2-1]
\end{tikzcd}\]

\end{lemma}
\begin{proof}
The first composition is given by 
\[\begin{tikzcd}
	&& X \\
	& {X\times_Z X} & {} & {X\times_Z X} \\
	X && {X\times_Z X \times_Z X} & {} & X.
	\arrow["{\Delta\times id_X}", from=2-2, to=3-3]
	\arrow["{id_X \times \Delta}"', from=2-4, to=3-3]
	\arrow["{id_X \times \pi}", from=2-4, to=3-5]
	\arrow["{\pi\times id_X}"', from=2-2, to=3-1]
	\arrow["\Delta"', from=1-3, to=2-2]
	\arrow["\Delta", from=1-3, to=2-4]
\end{tikzcd}\]
As $(\pi \times id_X)\circ \Delta=id_X= (id_X \times \pi)\circ \Delta$ this proves the claim. The claim for the second composition follows by the same argument.
\end{proof}

\begin{remark}\label{* and ! linearity}
    Let  $D : Corr(C,E) \to \PrL$ be a $6$-functor formalism, $Y,W, Z\in C$ with morphisms $\pi_Y: Y \to Z$ and $\pi_W: W \to Z$ and  $f: Y \to W$ a map in $E$ such that $\pi_Y =\pi_W \circ f$. Then the functors $f_!$ and $f^*$ are $D(Z)$-linear with respect to the $D(Z)$-module structures on $D(Y)$ and $D(W)$ induced by the functors $\pi_Y^*$ and $\pi_W^*$, respectively (cf. \Cite[Lecture 3, Remark 3.13]{6functors}).
\end{remark}

\begin{corollary}\label{D(complicated) equals id}
Let  $D : Corr(C,E) \to \PrL$ be a $6$-functor formalism, $\pi : X \to Z$ a map in $E$ and $\Delta : X \to X\times_Z X$ the diagonal map. Then we have isomorphisms of functors
\[
(id_X \times \pi)_!(id_X \times \Delta)^*(\Delta \times id_X)_!(\pi \times id_X)^* \cong id_{D(X)} \hspace{3mm}, \
(\pi \times id_X)_!(\Delta \times id_X)^*(id_X \times \Delta)_!(id_X \times \pi)^* \cong id_{D(X)}.
\]

\end{corollary}
\begin{proof}
Apply the functor $D(-)$ to the diagrams of \Cref{funnydiagrams} and use that they represent the identity morphisms in $\mathrm{Corr}(C,E)$.
\end{proof}

For  $\pi : X \to Z$ our fixed morphism, we introduce the notation $m \coloneq m_{X,X}$, $m_{12} \coloneq m_{X,X\times_Z X}$, $m_{21}\coloneq m_{X\times_Z X,X}$. As the next two propositions show, in order to deduce that $D(X)$ is dualisable as a $D(Z)$-module, it is sufficient for the map $m$ to be an equivalence.   
\begin{proposition}\label{magicmonsterdiagram proposition}
Let $D : Corr(C,E) \to \PrL$ be a $6$-functor formalism  and let $\pi: X \to Z$ be a morphism in $E$. If the functor 
 \[
 m: D(X)\otimes_{D(Z)}D(X) \to D(X\times_Z X)
 \]
 is an equivalence, then $D(X)$ is self-dual as a $D(Z)$-module.
\end{proposition}
 \begin{proof}
 We introduce the functors
 \begin{align*}
 t \coloneq & \pi_! \Delta^* : D(X\times_Z X) \to D(Z) \\
 u \coloneq & \Delta_! \pi^* : D(Z) \to D(X \times_Z X)
 \end{align*}
 which are $D(Z)$-linear by \Cref{* and ! linearity}. Using the $D(Z)$-linear equivalence  $m$, we can define the evaluation $ev$ and coevaluation $co$ by the diagrams 
\[\begin{tikzcd}
	{D(X\times_Z X)} & {D(Z)} & {D(X\times_Z X)} & {D(X)\otimes_{D(Z)}D(X)} \\
	{D(X)\otimes_{D(Z)}D(X)} && {D(Z)}
	\arrow["t", from=1-1, to=1-2]
	\arrow["m", from=2-1, to=1-1]
	\arrow["ev"', from=2-1, to=1-2]
	\arrow["u", from=2-3, to=1-3]
	\arrow["\sim"', from=1-4, to=1-3]
	\arrow["m", from=1-4, to=1-3]
	\arrow["co"', from=2-3, to=1-4]
\end{tikzcd}\]

\noindent  In order to see that $D(X)$ is dualizable as a $D(Z)$-module, we need to verify that 
 \[
 (ev\otimes id_{D(X)})\circ (id_{D(X)} \otimes co) \cong id_{D(X)}.
 \]
 To see this, consider the following diagram: 
\begin{figure*}[h]
 \tikzcdset{scale cd/.style={every label/.append style={scale=#1},
    cells={nodes={scale=#1}}}}

\begin{flushleft}
\[\begin{tikzcd}[scale cd=0.75]\label{bigdiagram}
	&& {D(X)\otimes_{D(Z)} D(X) \otimes_{D(Z)} D(X)} \\
	\\
	& {D(X)\otimes_{D(Z)}D(X\times_Z X)} && {D(X\times_Z X)\otimes_{D(Z)}D(X)} & {} \\
	{D(X)} && {D(X\times_Z X \times_Z X)} && {D(X)}
	\arrow["{(id \times \Delta)_!(id \times \pi)^*}"', from=4-1, to=4-3]
	\arrow["{(\pi \times id)_!(\Delta \times id)^*}"', from=4-3, to=4-5]
	\arrow["{m_{21}}"'{pos=0.6}, from=3-4, to=4-3]
	\arrow["{t\otimes id_{D(X)}}"'{pos=0.4}, from=3-4, to=4-5]
	\arrow["{m_{12}}"{pos=0.6}, from=3-2, to=4-3]
	\arrow["{id_{D(X)}\otimes u}"'{pos=0.6}, from=4-1, to=3-2]
	\arrow["{id_{D(X)}\otimes co}", curve={height=-30pt}, from=4-1, to=1-3]
	\arrow["{id_{D(X)} \otimes m}"{pos=0.3}, from=1-3, to=3-2]
	\arrow["{ev \otimes id_{D(X)}}", curve={height=-30pt}, from=1-3, to=4-5]
	\arrow["{m\otimes id_{D(X)}}"'{pos=0.3}, from=1-3, to=3-4]
\end{tikzcd}\]
\end{flushleft}
\caption{}
\end{figure*}

Note that we suppressed the equivalence $m_{X,Z}: D(X)\otimes_{D(Z)}D(Z) \cong D(X)$ in the lower left corner and the equivalence $m_{Z,X}: D(Z)\otimes_{D(Z)}D(X)\cong D(X)$ in the lower right corner of Figure 1 above.
First observe that the outer triangles involving $ev$ and $co$ commute by definition of $ev$ and $co$. Next we show the commutativity of the lower triangles, which follows by the naturality of $m$: Indeed, let $f: Y \to W$ be any morphism, then it follows from \Cref{naturality of m1} that the following diagram commutes

\begin{equation}\label{naturality of m 2}
\begin{tikzcd}
	{D(Y)\otimes_{D(Z)}D(X)} & {D(W)\otimes_{D(Z)}D(X)} \\
	{D(Y\times_Z X)} & {D(W\times_Z X).}
	\arrow["{f_! \otimes id_{D(X)}}", from=1-1, to=1-2]
	\arrow["{(f\times id_X)_!}", from=2-1, to=2-2]
	\arrow["{m_{Y,X}}"', from=1-1, to=2-1]
	\arrow["{m_{W,X}}", from=1-2, to=2-2]
\end{tikzcd}
\end{equation}

Combining \Cref{naturality of m1} and \Cref{naturality of m 2}, we see that each of the two lower triangles in Figure 1 commute. Finally, one easily verifies that the central square commutes by definition of the morphisms $m_{i,j}$. Thus we infer that 
\[
(ev\otimes id_{D(X)}) \circ (id_{D(X)} \otimes co) \cong (\pi \times id)_!(\Delta \times id)^*(id \times \Delta)_!(id \times \pi)^*.
\]
But this is isomorphic to the identity by \Cref{D(complicated) equals id}. The other identity 
\[
(id_{D(X)} \otimes ev)\circ (co \otimes id_{D(X)})\cong id_{D(X)}
\]
follows by an analogous argument. This shows that $D(X)\in \Mod_{D(Z)}(\PrL)$ is self-dual.
\end{proof}

\begin{corollary}\label{Künnethmorphismimplies dualisable}
    Let $f: X\to Z$ be Künneth. Then $D(X)\in \Mod_{D(Z)}(\PrL)$ is dualisable.
\end{corollary}
\begin{proof}
    Since $f$ is Künneth we have in particular that the morphism $m_{X,X}: D(X)\otimes_{D(Z)}D(X) \to D(X\times_Z X)$ is an equivalence. The result then follows  by \Cref{magicmonsterdiagram proposition}.
\end{proof}

\begin{corollary}\label{generalequivalence is given by a kernel}
Let $D : Corr(C,E) \to \PrL$ be a $6$-functor formalism  and let $X \to Z$ be a morphism in $C$. Assume that the functor $m: D(X)\otimes_{D(Z)}D(X) \to D(X\times_Z X)$ is an equivalence. Then we have an equivalence 
\[
\mathrm{Fun}_{D(Z)}^{L}(D(X),D(X))\cong D(X \times_Z X).
\]
 \end{corollary}
 \begin{proof}
We claim that there are equivalences
 \[
  \mathrm{Fun}_{D(Z)}^{L}(D(X),D(X))\cong  \mathrm{Fun}_{D(Z)}^{L}(D(Z), D(X)\otimes_{D(Z)} D(X))\cong D(X)\otimes_{D(Z)} D(X).
 \]
Indeed, by \Cref{magicmonsterdiagram proposition}, we know that $D(X)$ is self dual as a $D(Z)$ module which shows the first equivalence. The second equivalence is natural.
 \end{proof}

\subsection{Descent properties of Künneth morphisms}\label{Properties of Künneth morphisms}

In this subsection we will be interested in the interaction of the notion of Künneth morphisms with descent of the $6$-functor formalism $D: Corr(C,E) \to \PrL$. We denote by 
\[
D^*: C^{op} \to \PrL
\]
the functor given by $D^*(X)=D(X)$ for $X\in C$ and $D^*(f)=f^*$ for morphisms $f$ in $C$ and
\[
D^!: C_E^{op} \to \Cat
\]
the functor given by $D^!(X)=D(X)$ for $X\in C$ and $D^!(f)=f^!$ for morphisms $f$ in $E$ (c.f \Cite[Definition 3.1.4 and 3.2.3]{heyer20246}). In the following we use the notion of sieves as defined in \Cite[Definition A.4.1]{heyer20246}. 
\begin{definition}
 Let $C,D$ be categories and assume that all limits exist in $D$, let $\mathcal{U} \subset C_{/U}$ a sieve on an object $U\in C$ and $F: C^{op}\to D$ a functor. 
We say $F$ descends along $\mathcal{U}$ if the natural map 
     \[
     F(U) \to \lim_{V \in \mathcal{U}^{op}}F(V)
     \]
     is an isomorphism. We say $F$ descends universally along  $\mathcal{U}$ if it descends along every pullback of $\mathcal{U}$.
\end{definition}
Let us recall the definition of a universal $!$- and $*$-covers. 
\begin{definition}\Cite[Definition 3.4.6]{heyer20246}
 Let $D : Corr(C,E) \to \PrL$ be a  $6$-functor formalism. 
 \begin{enumerate}
     \item We say that a sieve $\mathcal{U} \subset (C_E)_{/U}$ is a small $!$-cover if it is generated by a small family of maps $(U_i \to U)_i$ and the functor $D^!$ descends along $\mathcal{U}$. We say $\mathcal{U}$ is a universal $!$-cover if for every map $V \to U$ in $C$ the family $(U_i\times_U V \to V)_i$ generates a small $!$-cover.
     \item We say that a sieve $\mathcal{U}$ in $C$ is a (universal) $*$-cover if $D^*$ descends (universally) along $\mathcal{U}$.
     \item The $D$-topology is the site on $C$ where covers are canonical covers which have universal $!$-and $*$-descent.
 \end{enumerate}
\end{definition}
For our purposes the following two classes of examples of $!$-covers will be particularly important.

\begin{proposition}\label{smooth and descendable D-covers}
Let $D : Corr(C,E) \to \Pr_{st}^L$ be a stable $6$-functor formalism, $f:X\to Y $ a morphism in $E$. Assume either of the following conditions are satisfied.
\begin{enumerate}
    \item (smooth cover) The object $1_X$ is $f$-smooth and the functor $f^*$ is conservative.
    \item (descendable cover) The object $1_X$ is $f$-proper and $f_*(1_X)\in \CAlg(D(Y))$ is descendable.
\end{enumerate}
Then the arrow $f : X \to Y$ is a universal $!$- and a universal $*$-cover.
\end{proposition}
\begin{proof}
This is \Cite[Proposition 6.18, Proposition 6.19] {6functors}.
\end{proof}
We will call a morphism $f:X\to Y $ in $E$ a smooth (respectively descendable) $!$-cover if it satisfies  condition $1$ (respectively condition $2$) of \Cref{smooth and descendable D-covers}.
As the following proposition shows, for $f: X \to Y$ a descendable or smooth $!$-cover, the category $D(Y)$ can be made more explicit.

\begin{proposition}\Cite[Proposition 3.1.27]{camargo2024analytic}\label{modulcomodul}
    Let $D : Corr(C,E) \to \Pr_{st}^L$ be a stable $6$-functor formalism and  $f: X \to Y$ be a morphism in $E$. 
 \begin{enumerate}   
  \item   If $f$ is a descendable $!$-cover we have 
    \[
    D(Y)\cong \coMod_{f^*f_*}D(X).
    \]
 \item    If $f$ is a smooth $!$-cover we have 
    \[
    D(Y)\cong \Mod_{f^!f_!}D(X).
    \]
\end{enumerate}
\end{proposition}

\begin{definition}\Cite[Definition A.4.5]{heyer20246}
$D : Corr(C,E) \to \PrL$ be a $6$-functor formalism and  $(f_i: {X}_i\to X)_{i\in I}$ be a universal $!$-cover of $X\in C$. We denote by $\Delta_I \to \Delta$ the right fibration associated with the functor $\Delta^{op}\to \Ani ,\  [n]\mapsto I^{n+1}$ by unstraightening. An object in $\Delta_I$ is a pair $([n]\in \Delta,i_\bullet\in I^{n+1})$ and a morphism $([n],i_\bullet)\to ([m],j_\bullet) $ in $\Delta_I$ is a morphism $\alpha : [n] \to [m]$ in $\Delta$ such that $i_k=j_{\alpha(k)}$ for all $k\in [n]$. We denote by $f_{\bullet} : X_{\bullet}\to X$ the Cech-nerve given by morphisms $f_{[n],i_\bullet} : X_{[n],i_\bullet}\coloneqq X_{i_0}\times_X ... \times_X X_{i_n} \to X $.
\end{definition}

We will now show one of our main results, which is that the property of a morphism being Künneth can by checked $!$-locally on the source. 
\begin{proposition}\label{dualisability is !-local on source}
Let $D : Corr(C,E) \to \PrL$ be a $6$-functor formalism and $g: X \to Z$ be any morphism in $E$. Let $(f_i: {X}_i\to X)_{i\in I}$ be a universal $!$-cover. If $g\circ f_{[n],i_\bullet} : X_{i_0}\times_X ... \times_X X_{i_n} \to Z$ is Künneth for all $([n],i_\bullet)\in \Delta_I$, then $g$ is Künneth.
\end{proposition}
\begin{proof}
By $!$-codescend, we have an equivalence $D(X)\cong \colim_{!,([n],i_\bullet)\in \Delta_I} D(X_{[n],i_\bullet})$ in $\PrL$. Using that the forgetful functor $\Mod_{D(Z)}(\Pr^L)\to \PrL$ is conservative and \Cite[Lemma 3.2.5]{heyer20246} (by passing to the induced geometric setup $(C_{/Z},E)$ of the slice category $C_{/Z}$) we see that this is also an equivalence in $\Mod_{D(Z)}(\PrL)$. Let $W \to Z$ be an arbitrary map. We claim that the following equivalences hold
\begin{align*}
    D(X)\otimes_{D(Z)} D(W)& \cong \colim_{!,([n],i_\bullet)\in \Delta_I} D(X_{[n],i_\bullet})\otimes_{D(Z)}D(W) \\
    & \cong \colim_{!,([n],i_\bullet)\in \Delta_I}  D(X_{[n],i_\bullet} \times_Z W) \\
    & \cong  D(X\times_Z W).
\end{align*}
In the first line we used that the functor $D(W) \otimes_{D(Z)} (-)$ commutes with colimits (see \Cref{tensor in Prl commutes with colimits}). The second line follows since the morphisms $X_i \to Z$ are Künneth. The last line follows by $!$-descent by noting that $(g_i :X_i\times_Z W \to X\times_Z W)_{i\in I}$ is a $!$-cover, as it is the base change of the universal $!$-cover $(X_i \to X)_{i\in I}$ and by observing that the morphisms $g_{[n],i_\bullet}: X_{[n],i_\bullet} \times_Z W \to X\times_Z W $ are given by the base change of $f_{[n],i_\bullet} : X_{[n],i_\bullet} \to X$ by $X\times_ZW \to X$. We now show that the equivalences 
\[
m_{{[n],i_\bullet}}: D(X_{[n],i_\bullet})\otimes_{D(Z)}D(W)\cong D(X_{[n],i_\bullet}\times_{Z}W)
\] 

are compatible with the colimit systems: By considering \Cref{naturality of m 2} for the morphims $h : X_{[n],i_\bullet} \to X_{[m],j_\bullet}$ in the colimit system we observe that the following diagram commutes.

\[\begin{tikzcd}
	{D(X_{[n],i_\bullet})\otimes_{D(Z)}D(W)} & {D(X_{[m],j_\bullet})\otimes_{D(Z)}D(W)} \\
	{D(X_{[n],i_\bullet}\times_Z W)} & {D(X_{[m],j_\bullet}\times_Z W).}
	\arrow["{h_{!}\otimes id_{D(W)}}", from=1-1, to=1-2]
	\arrow["{(h \times id_{W})_!}", from=2-1, to=2-2]
	\arrow["{m_{{[n],i_\bullet}}}"', from=1-1, to=2-1]
	\arrow["{m_{[m],j_\bullet}}", from=1-2, to=2-2]
\end{tikzcd}\]

This shows the compatibility of the equivalences $m_{{[n],i_\bullet}}$ with the the colimit systems. 
\end{proof}

Similarly, the notion of a morphism being Künneth can be checked $!$-locally on the target. 
\begin{proposition}
    
\label{dualisability is !-local on target}
 Let $D : Corr(C,E) \to \PrL$ be a $6$-functor formalism and $f: X \to Z$ be any morphism in $E$ and $(g_i: Z_i\to Z)_{i\in I}$ be a universal $!$-cover. Denote by  $f_{([n],i_\bullet)}: X_{([n],i_\bullet)}\coloneqq X \times_Z Z_{[n],i_\bullet} \to Z_{[n],i_\bullet}$ the base change of $f$ by $g_{[n],i_\bullet}: Z_{[n],i_\bullet} \to Z$. If $f_{[n],i_\bullet}$ is Künneth for all $([n],i_\bullet)\in \Delta_I$, then $f$ is Künneth.
\end{proposition}
\begin{proof}
Let $W \to Z$ be an arbitrary map and $W_{([n],i_\bullet)}\coloneqq W \times_Z Z_{[n],i_\bullet} \to Z_{[n],i_\bullet}$ the base change. Using $!$-descent for  the $!$-covers $X_{([n],i_\bullet)}\to X$,$Z_{([n],i_\bullet)}\to Z$, $W_{([n],i_\bullet)}\to W$ we obtain equivalences 
\begin{align*}
    D(X)\otimes_{D(Z)}D(W) \cong &\colim_{!,([n],i_\bullet)\in \Delta_I^{op}}D(X_{([n],i_\bullet)}) \otimes_{\colim_{!,([l],j_\bullet)\in \Delta_I^{op}}D(Z_{([l],j_\bullet)})}\colim_{!,([m],k_\bullet)\in \Delta_I^{op}}D(W_{([m],k_\bullet)}) \\
    \cong & \colim_{!,([n],i_\bullet)\in \Delta_I^{op}}D(X_{([n],i_\bullet)}\times_{Z_{([n],i_\bullet)}} W_{([n],i_\bullet)}) \\
    \cong & \colim_{!,([n],i_\bullet)\in \Delta_I^{op}}D(X\times_Z W \times_Z Z_{([n],i_\bullet)}) \\
    \cong & D(X\times_Z W)
\end{align*}
Here we used $!$-descent for the $!$-covers $X_{([n],i_\bullet)}\to X$,$Z_{([n],i_\bullet)}\to Z$, $W_{([n],i_\bullet)}\to W$ in the first line, the fact that  the category $\Delta_I^{op}$ is sifted in the second line, the fact that $f_{[n],i_\bullet}$ is Küneth for all $([n],i_\bullet)\in \Delta_I$ in the third line and  $!$-descent for $(X\times_Z W \times_Z Z_{([n],i_\bullet)} \to X\times_Z W$ in the last line.
\end{proof}
 In view of our applications to analytic stacks we investigate the following situation of an extension of a $6$-functor formalism.
\begin{proposition}\Cite[Theorem 3.4.11]{heyer20246}\label{E' properties}
    Let $D: Corr(C,E) \to \PrL$ be a sheafy $6$-functor formalism, where $C$ is a subcanoncial site. Then there is a collection of edges $\tilde{E}$ in $\mathcal{X}\coloneqq \Sh(C)$ with the following properties:
\begin{enumerate}[(i)]
\item The inclusion $C \hookrightarrow \mathcal{X}$ defines a morphism of geometric setups $(C,E) \to (\mathcal{X},\tilde{E})$ and $D$ extends uniquely to a $6$-functor formalism $\tilde{D}$ on $ (\mathcal{X},\tilde{E})$.
    \item $\tilde{E}$ is $*$-local on the target: Let $f: X \to Y$ in $C$ be map for which the pullback to any object in $C$ lies in $E$, then $f$ lies in $\tilde{E}$. 
    \item $\tilde{E}$ is $!$-local: Let $f : X \to Y$ be a map which is $!$-local on the source or target in $\tilde{E}$, then $f$ lies in $\tilde{E}$.
    \item $\tilde{E}$ is tame: Every map $f : X \to Y$ in $\tilde{E}$ with $Y\in C$ is $!$-locally on the source in $E$.
\end{enumerate}
Moreover, there is a minimal choice of $\tilde{E}$.
\end{proposition}

\begin{theorem}\label{Künnethextension-main theorem}
    Let $D: Corr(C,E) \to \PrL$ be a sheafy $6$-functor formalism on a subcanonical site $C$. Let $S\in C$ and assume that $D_S: ((C_E){{_{/S}}})^{op} \to \Mod_{D(S)}(\PrL)$ is symmetric monoidal. Then the extension $\tilde{D}_S: ((\Shv(C)_{\tilde{E}}){{_{/S}}})^{op}\to \Mod_{D(S)}(\PrL) $ is symmetric monoidal.
\end{theorem}
\begin{proof}
We review the inductive construction of $\tilde{E}$ in the proof of \Cite[Theorem 3.4.11.]{heyer20246} and show that the property of morphisms being Künneth is preserved under each extension step. Let $E_0$ be the class of morphisms in $\mathcal{X}$ which are representable in $E$ and by $A$ the collections of edges $E'$ such that $E'$ satisfies $(i)$ and $(iv)$ in \Cref{E' properties}. For $E' \in A$, let $(E')'\coloneqq E''$ be the class of morphisms in $\mathcal{X}$ which are $D^{!'}$-locally in $E'$. For $E'\in A$, we denote by $E_S'$ the class of morphisms $X \to S$ in $E'$ and by $Ku(E_S')\subset E_S'$ the class of morphisms $f: X \to S$ in $E'_S$ which are Künneth. Note that $E_0 \in A$. We first observe that $Ku(E_{0,S})=E_{0,S}$. Indeed, by representability we have that any $f : X \to S$ in $E_0$ lies in $E$ since $S\in C$, thus $f$ is Künneth since $D_S: ((C_E){{_{/S}}})^{op} \to \Mod_{D(S)}(\PrL)$ is symmetric monoidal. We show the following claim.
\[
(\ast) \ \text{If $E'\in A$ and $Ku(E_S')= E_S'$, then $Ku(E_S'')= E_S''$}
\]

Indeed, let $f: X \to S$ be an element in $E_S''$. By construction, there exists a small $!-$cover $(X_i\to X)_{i\in I}$ with $g_i:X_i\to X$ in $E'$ and the composition $f\circ g_i \in E'$. By stability under base change and composition of $E'$ we have that all the morphisms $f\circ g_{[n],i_\bullet}: X_{[n],i_\bullet} \to S$ are in $E'$. It thus follows that $f$ is Künneth by \Cref{dualisability is !-local on source} and thus $Ku(E_S'')=E_S''$. \\ For $E'\in A$, we define $E'_!\coloneqq \bigcup_{n\geq2}E^n$ where $E^2\coloneqq E''$ and $E^{n+1}\coloneqq (E^n)'$ for $n \geq 2$. By \Cite[Theorem 3.4.11]{heyer20246} we have that $E'_! \in A$ and satisfies $(iii)$. We have the following claim:

\[
(\ast') \ \text{If $E'\in A$ and $Ku(E_S')= E_S'$, then $Ku({E_{!S}'})= E'_{!S}$}
\]
Indeed, let $f: X \to S$ be in $E'_{!S}$. By definition, we know that $f\in E^n_S$ for some $n\ge2$, thus $f$ is Künneth by applying $(\ast)$ repeatedly. \\

For $E' \in A$ satisfying $(iii)$, we follow Heyer-Mann and define $E'_*$ to be the class of morphisms $X \to Y$ in $\mathcal{X}$ which after pullback under any morphism $U \to Y$ for $U\in C$ lie in $E'$. By \Cite[Theorem 3.4.11]{heyer20246} this defines again an element $E'_* \in A$. Since $E_{*S}'=E_S'$, we obtain the following claim:
\[
(\ast \ast) \ \text{If $E'\in A$  and $Ku(E_S')= E_S'$, then $Ku(E_{*S}')= E_{*S}'$}
\]
By construction we have $\tilde{E}\coloneqq \bigcup_{n \geq 1} \tilde{E}^n$, where $\tilde{E}^1\coloneqq E_{0!*}$ and $\tilde{E}^{n+1}\coloneqq (\tilde{E}^{n})_{!*}$. It follows from $(\ast),(\ast')$ and $(\ast \ast)$ that $Ku(\tilde{E}_S)=\tilde{E}_S$.
\end{proof}
\begin{remark}
\Cref{dualisability is !-local on source} together with \Cref{modulcomodul} can be used to construct dualisable comodule categories: For example, consider a $6$-functor formalism $D$ on $(\Sh(C),\tilde{E})$ as above and assume that the restriction of $D$ to $(C,E)$ satisfies Künneth. Let $G \in \Sh(C)$ be a commutative group object such that $h: G\to \ast$ is $D$-cohomologically proper and the natural map $f: \ast \to \ast/G$ is a descendable $D$-cover. We note that in this case $h_*\cong h_!$ satisfies the projection formula. Then using \Cref{modulcomodul} and the proper base change along the cartesian square 
\[\begin{tikzcd}
	G & \ast \\
	\ast & {\ast/G}
	\arrow["h", from=1-1, to=1-2]
	\arrow[from=1-1, to=2-1]
	\arrow["f", from=1-2, to=2-2]
	\arrow["f"', from=2-1, to=2-2]
\end{tikzcd}\]
we obtain an equivalence $D(\ast/G)\cong \coMod_{h_*(1)}D(\ast)$. In fact, this equivalence is stable under arbitry base change, that is for any $\pi: Y \to \ast$ we obtain an equivalence 
\[
D(\ast/G \times_\ast Y)\cong \coMod_{\pi^*h_*(1)}D(Y).
\]
Using  \Cref{dualisability is !-local on source} we note that we obtain a Künneth formula for comodules: 
\[
\coMod_{h_*(1)}D(\ast) \otimes_{D(\ast)} D(Y)\cong \coMod_{\pi^*h_*(1)}D(Y). 
\]
Although the corresponding statement for smooth $D$-covers reduces to a well-known statement about categories of modules (cf. \Cite[Proposition 4.1]{ben2010integral}) a corresponding abstract statement for comodules is unknown to the author.
\end{remark}

\begin{theorem}\label{Künnethextension-main corollary}
    Let $D: Corr(C,E) \to \PrL$ be a $6$-functor formalism which is Künneth and consider the extension $\tilde{D}: Corr(\Shv(C),\tilde{E})\to \PrL$. Assume that $S\in \Shv(C)$ admits a $!$-cover $(S_i\to S)_{i\in I}$ with $S_i \to S$ in $E_0$ and $S_i \in C$ for all $i\in I$. Then any morphism $f : X \to S$ in $\tilde{E}$ is Künneth.
\end{theorem}
\begin{proof}
Since $S_i \to S$ in $E_0$ we have that $S_{[n],i_\bullet}\in C$ for all $([n],i_\bullet) \in \Delta_I$. The result follows now by \Cref{Künnethextension-main theorem} and \Cref{dualisability is !-local on target}. 
\end{proof}
\begin{remark}
   The previous corollary applies for example if $S \in \Shv(C)_{\tilde{E}/\ast}$ has representable diagonal, $\Delta_S : S \to S \times S \in E_0$. In the context of derived algebraic geometry and with $D(-)$ given by $\QCoh(-)$ this  corresponds to the affineness condition on the diagonal of the base in the definition of perfect stacks (c.f \Cite[Definition. 3.2]{ben2010integral}). One could slightly generalise this to the affineness condition to allow "quasi-affine" diagonals (to reproduce the results of \Cite[Theorem 1.0,3]{stefanich2023tannaka}) using the fact that open immersion are Künneth \Cref{open immersion and closed immersions are Künneth}.
\end{remark}
We include the following lemma
\begin{lemma}\label{fibreproductcommutes with colimits}
    Let $\Sh(C)$ be the category of sheaves of anima on a small site $C$, $Z$ an object in $\Sh(C)$ and $X,Y\in \Sh(C)_{/Z}$ . Then the functor given by taking fibre products
\[
X\times_Z (-) : \Sh(C)_{/Z} \to \Sh(C)_{/Z} ,  \ Y \mapsto X\times_Z Y
\]
commutes with all colimits. In particular, it has a left adjoint which is given by the Hom-stack  
\[
\Hom_Z(X,Y) : C_{/Z} \to \Ani ,  \ U \mapsto \text{Hom}_Z(X \times_Z U,Y).
\]
Here, the category $C_{/Z} \hookrightarrow \Sh(C)_{/Z}$ is given by pullback of the forgetful functor $\Sh(C)_{/Z}\to \Sh(C)$ along the fully faithful inclusion $C \hookrightarrow \Sh(C)$ via Yoneda.
\end{lemma}
\begin{proof}
The fact that the functor $X\times_Z (-)$ commutes with all colimits is immediate as $\Sh(C)$ and hence $\Sh(C)_{/Z}$ are $\infty-$topoi. Since $\Sh(C)_{/Z}$ is presentable, we have a right adjoint $L_X$ which is given on  $U\in C_{/Z}$ by
\begin{align*}
L_X(U)\cong & \text{Hom}_Z(U, L_X) \\
\cong & \text{Hom}_Z(X \times_Z U,Y) \\
\cong & \Hom_Z(X,Y)(U).
\end{align*}
\end{proof}

\begin{lemma}\label{Künneth remains künneth under extension}
  Let $D_0 : Corr(C_0,E_0) \to \PrL$ be a $6$-functor formalism  on a subcanonical site $C_0$ and consider an extension $D : Corr(\Sh(C_0), \tilde{E}_0) \to \PrL$ as in \Cref{E' properties}. Then $f : X \to Y$ in $E_0$ is Künneth with respect to $D_0$ if and only if $f$ is Künneth with respect to $D$.
\end{lemma}
\begin{proof}
If $f$ is Künneth with respect to $D$ then it is also Künneth with respect to $D_0$ since $C_0 \subset \Sh(C_0)$. Conversely, since $f\in E_0$ we obtain $f\in \tilde{E}_0$ by the definition of the extension $\tilde{E}_0$. Let $W \to Z$ be a morphism in $\Sh(C_0)$ and write $D(W) \cong \lim_{n\in J} D_0(W_n)$ with $W_n\in C_0$. By \Cref{Künnethmorphismimplies dualisable} we obtain that $D_0(X) \in \Mod_{D_0(Z)}(\PrL)$ is dualisable.  We claim that there is a chain of equivalences
\begin{align*}
D_0(X) \otimes_{D_0(Y)}D(W) & \cong \lim_{n\in J}D_0(X) \otimes_{D_0(Y)}D_0(W_n) \\
& \cong \lim_{n\in J} D_0(X\times_Y W_n) \\
& \cong D(X \times_Y W).
\end{align*}
Indeed, the first line follows by \Cref{tensor with dualisable object commmutes lim}, the second equivalence follows from the fact that $f$ is Künneth with respect to $D_0$ and the last one follows from \Cref{fibreproductcommutes with colimits}.
\end{proof}

\begin{remark}\label{remark diagobal and rightcancellativenenss}
Let $D : Corr(C,E) \to \PrL$ be a $6$-functor formalism. By \Cref{dualisability is stable under product},\Cref{dualisable stable under bc} and \Cref{dualisability is !-local on source} we see that the class $K$ of Künneth morphisms in $C$ has good properties such as closure under composition, arbitrary base change and $!$-locality on the source. Note however, that $(C,K)$ does not define a geometric setup in the sense of \Cite[Definition 2.1.1]{heyer20246}. The problem is that the class $K$ is not right cancellative in general, or in view of \Cite[Remark 2.1.2]{heyer20246} that for $f: X \to Y$ Künneth the diagonal $\Delta_f: X \to X\times_Y X$ is not necessarily Künneth. 
\end{remark}

\subsection{Consequences: Tannakian reconstruction}\label{Consequences: Tannakian lifting}
Let $D: Corr(C,E)\to \CAlg(\PrL)$ be a $6$ functor formalism and let $S \in C$. It is an interesting question whether one can recover maps $\text{Hom}_{C_{/S}}(X,Y)$ in $C_{/S}$ from the category of morphisms $\Fun_{D(S)}^{L,\otimes}(D(Y),D(X))$ in $\CAlg(\Mod_{D(S)}(\PrL))$. We introduce the following terminology.
\begin{definition}\label{tannakian morphism}
  Let $D: Corr(C,E)\to \CAlg(\PrL)$ be a $6$-functor formalism and $f: X \to S$ in $C$. We call $f$ Tannakian if the functor 
  \[
 \text{Hom}_{C_{/S}}(Y,X) \to \Fun_{D(S)}^{L,\otimes}(D(X),D(Y))
  \]
  is an equivalence for all $Y \in C_{/S}$. We call $D$ Tannakian if all morphisms in $C$ are Tannakian.
\end{definition}
Given a $6$-functor formalism $D: Corr(C,E)\to \CAlg(\PrL)$ on a subcanonical site $C$ which is Tannakian and Künneth, we will investigate which morphisms $f :X \to S$ in the extension $\tilde{D}: Corr(\Shv(C),\tilde{E})\to \CAlg(\PrL)$ given by \Cref{E' properties} are Tannakian.  Our arguments basically follow \Cite[Section 4]{stefanich2023tannaka} by showing that the class of $!$-covers in \Cref{Künnethextension-main corollary} satisfy $2$-descent. More generally we show the following lemma.
\begin{lemma}\label{shrik descent implies 2 descent}
Let $D: Corr(C,E)\to \CAlg(\PrL)$ be a $6$-functor formalism and $f: X' \to X$ in $E$ satisfying Künneth. Then 1) implies 2) implies 3):
\begin{enumerate}
    \item $f$ is a universal $!$-cover
    \item For all $M\in \Mod_{D(X)}\PrL$ we have that the natural map $M \cong \underset{[n]\in \Delta}{\lim}^* M \otimes_{D(X)}D(X'_{[n]})$ is an isomorphism. In particular $f$ satisfies $^*$-descent.
    \item $f$ satisfies $2$-descent, i.e. we have an isomorphism $\Mod_{D(X)}\PrL \cong \lim_{[n]\in \Delta} \Mod_{D(X'_{[n])})}\PrL $.
\end{enumerate}
\end{lemma}
\begin{proof}
We start by showing that 1) implies 2). Let $M\in \Mod_{D(X)}\PrL$, then we have equivalences 
\begin{align*}
M\cong & \Fun_{D(X)}^L(D(X),M) \\
\cong & \Fun_{D(X)}^L(\colim_{[n]\in \Delta^{op}}D(X'_{[n])}),M)\\
\cong & \lim_{[n]\in \Delta}\Fun_{D(X)}^L(D(X'_{[n])}),M)\\ 
\cong & \lim_{[n]\in \Delta}\Fun_{D(X)}^L(D(X), D(X'_{[n])})\otimes_{D(X)}M)\cong \lim_{[n]\in \Delta} M \otimes_{D(X)}D(X'_{[n])})
\end{align*}
Here we used that $f: X' \to X$ and hence also the basechanges $f: X'_{[n]} \to X$ are Künneth, which implies that $D(X'_{[n]})\in \Mod_{D(X)}\PrL$ are self-dual by  \Cref{Künnethmorphismimplies dualisable} in the fourth line. Taking $M=D(X)$ it follows that $f$ satisfies $*$-descent.\\
To see that 2) implies 3), note that the natural functor 
\[
B: \Mod_{D(X)}\PrL \to \lim_{[n]\in \Delta} \Mod_{D(X'_{[n]})}\PrL  \ \ , M \mapsto M_n\coloneqq M \otimes_{D(X)}D(X'_{[n]})
\]
admits a right adjoint given by $\lim_{[n]\in \Delta} \Mod_{D(X'_{[n]})}\PrL \to \Mod_{D(X)}\PrL$, $M_n \mapsto \lim_{[n]\in \Delta} f_{n,*}M_n $ where $f_{n,*}: \Mod_{D(X'_{[n]})}\PrL \to \Mod_{D(X)}\PrL $ is the forgetful functor. By 2), the unit of this adjunction is an equivalence. Thus, it suffices to verify that the counit is an equivalence. Since  $D(X'_{[n]})\in \Mod_{D(X)}\PrL$ are dualisable, it suffices to show that the canonical morphism 
\[
(\lim_{[n]\in \Delta} f_{n,*}M_n) \otimes_{D(X)}D(X'_{[m]}) \to  M_m
\]
is an equivalence. We calculate 
\begin{align*}
  (\lim_{[n]\in \Delta} f_{n,*}M_n) \otimes_{D(X)}D(X'_{[m]}) \cong  \lim_{[n]\in \Delta} (f_{n,*}M_n \otimes_{D(X)}D(X'_{[m]})) \cong M_m. 
\end{align*}
\end{proof}

We can now formulate an extension lemma for Tannkian morphisms. The idea of the proof is adapted from \Cite[Theorem 4.2.1]{stefanich2023tannaka}.
\begin{theorem}\label{tannakian lift}
    Let $C$ be a subcanonical site and $D: Corr(C,E)\to \CAlg(\PrL)$ be a Tannakian sheafy $6$-functor formalism satisfying Künneth. Consider its extension $D: Corr(\Shv(C),\tilde{E})\to \CAlg(\PrL)$. Let $f: X \to S$ be a morphism in $\tilde{E}$ and assume that there is a $!$-cover $g: S'\to S$ with $S'\in C$ and $g\in E_0$. Then $f$ is Tannakian.
\end{theorem}
\begin{proof}
  By iterating the following argument we can assume $S= \ast \in C$ the final object.  Consider the following commutative diagram 
\[\begin{tikzcd}
	{(C_{/X})^{op}} & {\CAlg(\Mod_{D(\ast)}(\PrL))_{D(X)/}} \\
	{C^{op}} & {\CAlg(\Mod_{D(\ast)}(\PrL))}
	\arrow["D", from=1-1, to=1-2]
	\arrow[from=1-1, to=2-1]
	\arrow[from=1-2, to=2-2]
	\arrow["D", hook, from=2-1, to=2-2]
\end{tikzcd}\]
Note that the lower horizontal functor is fully faithful since  $D: Corr(C,E)\to \CAlg(\Mod_{D(\ast)}(\PrL))$ is Tannakian by assumption. It suffices to show that this is a pullback diagram, the morphism 
\[
\text{Hom}_{C/\ast}(Y,X) \to \Fun_{D(\ast)}^{L,\otimes}(D(X),D(Y))
\]
arises as the fibre over $D(X) \in \CAlg(\Mod_{D(\ast)}(\PrL))$ of the left fibration $\CAlg(\Mod_{D(\ast)}(\PrL))_{D(X)/}\to \CAlg(\Mod_{D(\ast)}(\PrL))$.
Let $P$ be the pullback of the diagram, which we can identify as as full subcategory $\Mod_{D(X)}^{\text{aff}}\PrL$ of $\Mod_{D(X)}\PrL$ consisting of those $M\in \Mod_{D(X)}\PrL $ for which $M\cong D(Z)$ for some $Z\in C$. Let 
\[
\Theta : (C_{/X})^{op} \to  \Mod_{D(X)}^{\text{aff}}\PrL
\]
be the induced functor. We want to show it is an equivalence. Note that since $S' \to S'$ satisfies $!$-descent and is Künneth by \Cref{Künnethextension-main corollary}, this is also true for $X'\coloneqq X\times_S S' \to X$ by \Cref{dualisable stable under bc}. By \Cref{shrik descent implies 2 descent} we have $2$-descent for $X' \to X$ and since $g\in E_0$ we see that the equivalence from 3) in \Cref{shrik descent implies 2 descent} descents to an equivalence 
\[
\Mod_{D(X)}^{\text{aff}}\PrL \cong \lim_{[n]\in \Delta} \Mod_{D(X'_{[n]})}^{\text{aff}}\PrL.
\].
We thus obtain a diagram 
\[\begin{tikzcd}
	{(C_{/X})^{op}} && {\CAlg(\Mod_{D(X)}^{\text{aff}}\PrL)} \\
	{\underset{[n]\in \Delta}{\lim}(C_{/X'_{[n]}})^{op}} && {\underset{[n]\in \Delta}{\lim}\CAlg(\Mod_{D(X'_{[n]})}^{\text{aff}}\PrL)}
	\arrow["\Theta", from=1-1, to=1-3]
	\arrow[hook', from=1-1, to=2-1]
	\arrow[hook', from=1-3, to=2-3]
	\arrow["{\underset{[n]\in \Delta}{\lim}\Theta_n}"', from=2-1, to=2-3]
\end{tikzcd}\]
Note that for any $n$, the functor $\Theta_n$ is an equivalence, hence the limit $\lim_n\Theta_n$ is an equivalence. It follows that $\Theta$ is fully faithful. We now show essential surjectivity of $\Theta$. Let $D(X) \to D(Z)\in \CAlg(\Mod_{D(X)}^{\text{aff}}\PrL)$, then by fully faithfullness of the right vertical arrow and the equivalence $\underset{[n]\in \Delta}{\lim}\Theta_n$ we obtain compatible arrows $Z\times_X X'_{[n]} \to  X'_{[n]}$. But this gives an arrow $Z \to X$ by $!$-descent. so $\Theta$ is an equivalence.

\end{proof}

\newpage

\section{Künneth morphisms for analytic stacks}\label{$6$-functor formalism for analytic stacks}
\subsection{The $6$-functor formalism for analytic stacks}\label{6functorsforanalyticstacks}
In this subsection we recall the construction of the $6$-functor formalism of analytic stacks, following Clausen-Scholze \Cite[Lecture 16-20]{AnSt}. Since we are mainly interested in formal questions of $6$-functor formalisms, we will not recall all the relevant definitions and refer to \Cite{AnSt}, \Cite{clausen2022condensed}.   We denote by  $D(\mathbb{Z})\coloneqq D(\text{CondAb})$ be the derived category of light condensed abelian groups. 
\begin{definition}\label{defanring}\Cite[Lecture 8]{AnSt}
    An analytic ring $A=(A^{\tri},D(A))$ is the datum of a condensed animated ring $A^{\tri}$ and a fully faithful subcategory $D(A) \overset{j_A}{\subset} \Mod_{A^{\tri}}D(\mathbb{Z})$ such that
    \begin{enumerate}
       \item The category $D(A)$ is stable under all limits and colimits in $\Mod_{A^{\tri}}D(\mathbb{Z})$ and there is a left adjoint of the inclusion $j_A: D(A) \subset \Mod_{A^{\tri}}D(\mathbb{Z})$ which we denote by $(-)_A^\wedge$.
       \item For all $M \in D(\mathbb{Z})$ and $N \in D(A)$, we have that $R\underline{\text{Hom}}_{\mathbb{Z}}(M,N) \in D(A)$.
       \item The left adjoint $(-)_A^\wedge$ sends connective to connective objects.
       \item The underlying condensed animated ring $A^{\tri}$ lies in  $D(A)$.
    \end{enumerate}
\end{definition}

We denote by  $\AnR$ be the category of analytic rings, by $\Aff\coloneqq \AnR^{op} $ its opposite category and by $\SpecAn(-)\coloneqq (-)^{op}: \AnR \to\Aff$ the opposite functor. Condition (1) in \Cref{defanring} implies that $D(A)$ is an accessible localisation of $\Mod_{A^{\tri}}D(\mathbb{Z})$ and thus a presentable category. For any morphism $f: \SpecAn(B) \to \SpecAn(A)$ in $\Aff$, there is an associated  functor $f^* : D(A) \to D(B)$ in $\PrL$ which is defined by the diagram 
\[\begin{tikzcd}
	{\Mod_{A^{\tri}}D(\mathbb{Z})} & {\Mod_{B^{\tri}}D(\mathbb{Z})} \\
	{D(A)} & {D(B)}
	\arrow["{f^{\tri *}}", from=1-1, to=1-2]
	\arrow["{(-)_B^{\wedge}}", from=1-2, to=2-2]
	\arrow["{j_A}", hook, from=2-1, to=1-1]
	\arrow["{f^*}"', from=2-1, to=2-2]
\end{tikzcd}\]
Here we denote by $f^{\tri,*}\coloneqq B^{\tri} \otimes_{A^{\tri}}(-)$ the base change functor. Denoting the right adjoint of $f^*$  by $f_*$ we thus obtain a functor 
\begin{align*}
  D : \Aff^{op} \to \CAlg(\PrL), \ \SpecAn(A) &\mapsto D(A)   \\
   (f: \SpecAn(B) \to \SpecAn(A)) &\mapsto f^* : D(A) \to D(B)
\end{align*}
\begin{remark}
In order to apply \Cref{tannakian lift}, it will be convenient to work with  (condensed) $\mathbb{E}_{\infty}$-rings as opposed to (condensed) animated rings and to consider analytic $\mathbb{E}_{\infty}$-rings. We denote by $\text{AnCR}$ the category of analytic complete $\mathbb{E}_{\infty}$-rings as in \Cite[Definition 2.1.1]{camargo2024analytic}. Indeed, by \Cite[Lemma 2.1.3]{camargo2024analytic} the resulting $6$-functor formalism on $\text{AnCR}$ is Tannakian. Everything we say in the following will be correct as well for  $\text{AnCR}$ replacing $\AnR$.
\end{remark}
Our aim is now to apply the construction result \Cite[Proposition 3.3.3]{heyer20246}. As a first step we thus need an appropriate class of morphisms $I,P$ and $E$ in $\Aff$.  
\begin{definition}\Cite[lecture 17]{AnSt}
Let $\pi: X\to Y$ be a map in $\Aff$.
\begin{enumerate}
    \item We say $\pi$ is proper if $\pi_*$ commutes with small colimits and satisfies the projection formula.
    \item We say $\pi$ is an open immersion if it is a $D$-open immersion in the sense of \Cref{def open immersion for D Cop}.
    \item We say $\pi$ is $!$-able if it can by written as a composition $\pi=p \circ j$ of an open immersion $j$ and a proper map $p$.
\end{enumerate}
We denote by $E$ the class of $!$-able maps and by $I,P \subset E$ the class of open immersions and of proper maps, respectively.
\end{definition}
 We will need the following alternative description of proper maps. 
\begin{lemma}\label{alternative char. of open and proper in aff}
Let $\pi: \SpecAn(A) \to \SpecAn(B)$ be a map in $\Aff$. 
\begin{enumerate}
    \item $\pi$ is proper if and only if $A$ carries the induced analytic ring structure of $B$, i.e.\  $A=A^{\tri}_{/B}\coloneqq (A^{\tri},\Mod_{A^{\tri}}D(B))$. 
    \item If $\pi : X \to Y$ is proper then it satisfies proper base change, i.e. for any  Cartesian square
    \[\begin{tikzcd}
	{X'} & Z \\
	X & Y
	\arrow["{\pi'}", from=1-1, to=1-2]
	\arrow["{f'}", from=1-1, to=2-1]
	\arrow["f", from=1-2, to=2-2]
	\arrow["\pi", from=2-1, to=2-2]
\end{tikzcd}\]
the natural morphism of functors $f^*\pi_* \to f'^*\pi'_*$ is an equivalence. The base change of a proper map is again proper. Any composition of proper maps is proper.
   
\end{enumerate}
\begin{proof}
1) Assume $A$ carries the induced analytic ring structure. Then $\pi_* : \Mod_{A^{\tri}}D(B) \to D(B)$ is the forgetful functor, and hence commutes with colimits and is $D(B)$-linear and thus satisfies the projection formula.  Now assume that $\pi$ is proper. By the projection formula we have for any $N \in D(B)$ equivalences 
\begin{align*}
\pi_*(A^{\tri}\otimes_A \pi^*(N)) &\cong \pi_*(A^{\tri})\otimes_B N \\
& \cong  (\pi_*(A^{\tri}) \otimes_{B^{\tri}}N)^{\wedge}_B.
\end{align*}
However, the left hand side identifies to $\pi_*( (A^{\tri} \otimes_B N)^{\wedge}_A)$. We thus obtain that the element $A^{\tri}\otimes_B N \in \Mod_{A^{\tri}}D(B) $ lies automatically in $D(A)$, since the former is generated by colimits of the elements $A^{\tri}\otimes_B N$ for $N\in D(B)$  we obtain the desired equivalence $D(A) \cong \Mod_{A^{\tri}}D(B)$. \\
2)  Assume  that $\pi$ is proper. Let $f: \SpecAn(C^{\tri},D(C))\to \SpecAn(B^{\tri},D(B))$ any morphism of affine analytic stacks, then we obtain by the description of pushouts in $\AnR$  !!!REF!!! that the fibre product is given by $\SpecAn(A)\times_{\SpecAn(C)}\SpecAn(B)\cong \SpecAn((A^{\tri}\otimes_{B^{\tri}}C^{\tri},\Mod_{A^{\tri}\otimes_{B^{\tri}}C^{\tri}}D(C))$. Let $\tilde{\pi}$ be the base change of $\pi$ along $f$. The base change diagram reduces to the diagram 

\[\begin{tikzcd}
	{\Mod_{A^{\tri}\otimes_{B^{\tri}}C^{\tri}}D(C)} & {D(C)} \\
	{\Mod_{A^{\tri}}D(B)} & {D(B)}
	\arrow["{\tilde{\pi}_*}", from=1-1, to=1-2]
	\arrow["{\tilde{f}^*}"', from=2-1, to=1-1]
	\arrow["{\pi_*}", from=2-1, to=2-2]
	\arrow["{f^*}"', from=2-2, to=1-2]
\end{tikzcd}\]
which indeed commutes. In particular we see that the fibre product carries the induced analytic ring structure from $(C^{\tri},D(C))$, which shows that $\tilde{\pi}$ is proper by part 1).  The claim that a composition of proper maps is again proper can be checked easily from the definition of proper maps. 
\end{proof}

\end{lemma}

\begin{lemma}\label{pbc and stability under pullback of open imersion in aff}
Let $\pi: \SpecAn(A) \to \SpecAn(B)$ be a map in $\Aff$.  If $\pi$ is an open immersion then it satisfies proper base change. The base change of an open immersion is again an open immersion. Any composition of open immersions is an open immersion.
\end{lemma}
\begin{proof}
  Let $I\in D(B)$ be the idempotent algebra corresponding to the open immersion $\pi$, so $D(A) \subset D(B)$ identifies with the full subcategory of elements $M \in D(B)$ for which $I \otimes M\cong 0$. Let $f: \SpecAn(C) \to \SpecAn(B)$ be any map in $\Aff$. By the definition of pushouts in $\AnR$, we know that $C\otimes_B A$ is given by the completion of the analytic ring structure on $C^{\tri}$ such that  module $M \in \Mod_{C^{\tri}\otimes_{B^{\tri}} A^{\tri}}D(\mathbb{Z})$ lies in $D(C\otimes_B A)$ if and only its restriction to $\Mod_{C^{\tri}}D(\mathbb{Z})$ respectively $\Mod_{A^{\tri}}D(\mathbb{Z})$ lies in $D(C)$ respectively $D(A)$. Thus $M \in D(C)$ lies in $D(C \otimes_B A)$ if and only $M\otimes I\cong 0$. This shows that the base change is again an open immersion. The claim that proper base change holds follows by the construction the equivalence $D(C \otimes_B A)\cong \{M\in D(C)\mid f^*(I)\otimes M\cong 0 \}$. The claim that a composition of open immersions is again an open immersion can be checked easily from the definition of open immersions. 
\end{proof}

\begin{lemma}\label{AffE is geometric setup}
The pair $(\Aff, E)$ is a geometric setup
\end{lemma}
\begin{proof}
 Since open immersions and proper maps are stable under arbitrary base change in $\Aff$ (cf. \Cref{alternative char. of open and proper in aff}, \Cref{pbc and stability under pullback of open imersion in aff}), the same is true for any $!$-able map. We now check stability under composition. Let $f : \SpecAn(A) \to \SpecAn(B)$ and $g: \SpecAn(B) \to \SpecAn(C)$ be $!$-able maps and write $f=p_1 \circ j_1$  and $g= p_2 \circ j_2$ for $j_i \in I$ and $p_i\in P$. We consider the factorisation 
\[\begin{tikzcd}
	& {\SpecAn(A^{\tri}_{/C})} \\
	{\SpecAn(A)} && {\SpecAn(C)}
	\arrow["p", from=1-2, to=2-3]
	\arrow["\iota", from=2-1, to=1-2]
	\arrow["{g\circ f}"', from=2-1, to=2-3]
\end{tikzcd}\]
Here the map $p$ is clearly proper, since it is given by the induced analytic ring structure from $C$. It thus suffices to check that $\iota $ is an open immersion. Note that we have a cartesian diagram:
\[\begin{tikzcd}
	{\SpecAn(A^{\tri}_{/B})} & {\SpecAn(A^{\tri}_{/C})} \\
	{\SpecAn(B)} & {\SpecAn(B^{\tri}_{/C})}
	\arrow["{\tilde{j_2}}", hook, from=1-1, to=1-2]
	\arrow["{\tilde{p_1}}", from=1-1, to=2-1]
	\arrow["p_1", from=1-2, to=2-2]
	\arrow["j_2", hook, from=2-1, to=2-2]
\end{tikzcd}\]
 It thus follows that $\iota \cong {\tilde{j_2}} \circ j_1$ is an open immersion since it is a composition of open immersions (cf. \Cref{pbc and stability under pullback of open imersion in aff}).
\end{proof}

\begin{lemma}\label{suitable decomposition I,P}
The classes $I,P \subset E $ satisfy the following properties: 
\begin{enumerate}
    \item $I$ and $P$ contain all identity morphisms and are stable under compositions and arbitrary base change in $\Aff$.
    \item Every $f\in E$ can be written as a composition $f= p \circ j$ with $p\in P$, $j\in I$.
    \item  $I$ and $P$ are right cancellative. i.e. if $p, p\circ q \in I$ (respectively in $P$) then $q \in I$ (respectively in $P$)
    \item Every morphism $f \in I\cap P$ is $n$-truncated for some $n \geq -2$ (possibly depending on $f$).
\end{enumerate}
\end{lemma}
\begin{proof}
1) An isomorphism $X\cong Y$ in $\Aff$ is clearly an open immersion and proper.  The claim that $I$ and $P$ are stable under base change and composition is \Cref{alternative char. of open and proper in aff} and \Cref{pbc and stability under pullback of open imersion in aff}. The condition 2) is satisfied by definition of the class $E$. \\
3) Consider the diagram 
\[\begin{tikzcd}
	{\SpecAn(A)} & {\SpecAn(B)} \\
	& {\SpecAn(C)}
	\arrow["q", from=1-1, to=1-2]
	\arrow["{p\circ q}"', from=1-1, to=2-2]
	\arrow["p", from=1-2, to=2-2]
\end{tikzcd}\]

and consider the proper case. Note that for any proper map $p: \SpecAn(A^{\tri},\Mod_{A^{\tri}}D(B))\to \SpecAn((B^{\tri},D(B)))$ the functor $p_* : \Mod_{A^{\tri}}D(B) \to D(B)$ is the forgetful functor and hence conservative (cf. \Cref{alternative char. of open and proper in aff}). Assume $p, p\circ q$ are proper. Then $q_*$ commutes with all colimits since $p_*$ commutes with colimits and $p_*\circ q_*$ is conservative. The fact that $q_*$ satisfies the projection formula follows from \Cref{examples of internal leftadjoint}.\\
We now consider the case of open immersions. Let $j\coloneqq p\circ q$ and consider $D\coloneqq \mathrm{cofib}(j_!(1)\to 1)\in D(C)$ the idempotent algebra corresponding to the open immersion $j$. Then $p^*(D)\in D(B)$ is again idempotent and we claim that it defines an open immersion via \Cref{alternative char. of open and proper in aff}. We thus need to show 
\[
\mathrm{ker}(q^*)\cong \Mod_{p^*(D)} D(B)
\]
as full subcategories of $D(B)$. If $M\in \Mod_{p^*(D)} D(B)$, then $p^*(D)\otimes_B M\cong M$ so applying $q^*(-)$ we obtain $q^*p^*(D)\otimes_A q^*(M)\cong q^*(M)\cong 0$ since $q^*p^*(D)\cong j^*(D)\cong 0$, so $M \in \mathrm{ker}(q^*)$. Conversely, let $M\in \mathrm{ker}(q^*)$, and let $X\coloneqq \SpecAn(A), Y\coloneqq \SpecAn(B), Z\coloneqq \SpecAn(C) $. Consider the following cartesian diagram:

\[\begin{tikzcd}
	{X\times_ZY} & {Y\times_ZY} & Y \\
	X & Y & Z
	\arrow["{\tilde{q}}", from=1-1, to=1-2]
	\arrow["{\tilde{j}}", curve={height=-18pt}, from=1-1, to=1-3]
	\arrow["h", from=1-1, to=2-1]
	\arrow["\lrcorner"{anchor=center, pos=0.125}, draw=none, from=1-1, to=2-2]
	\arrow["{\tilde{p}}", from=1-2, to=1-3]
	\arrow["{\tilde{p}}", from=1-2, to=2-2]
	\arrow["\lrcorner"{anchor=center, pos=0.125}, draw=none, from=1-2, to=2-3]
	\arrow["p", from=1-3, to=2-3]
	\arrow["q", from=2-1, to=2-2]
	\arrow["j", curve={height=18pt}, from=2-1, to=2-3]
	\arrow["p", from=2-2, to=2-3]
\end{tikzcd}\]
We want to verify that $p^*(D)\otimes_B M\cong M$. Note that 
\begin{align*}
 p^*(D)\otimes_B M & \cong \mathrm{cofib}(p^*(j_!(1))\otimes_BM\to M) \\
 & \cong \mathrm{cofib}(\tilde{j}_!\tilde{j}^*(M)\to M).
\end{align*}
Here we used the symmetric monoidality of $p^*(-)$ in the first line and the fact that $\tilde{j}$ an open immersion, as it is a base change of the open immersion $j$ and hence satisfies the projection formula for $\tilde{j}_!$. It thus suffices to show that $\tilde{j}_!(\tilde{j}^*(M))\cong 0$, which follows by the following diagram chase: 
\begin{align*}
 \tilde{j}_!\tilde{j}^*(M)&\cong q_!h_!\tilde{q}^*\tilde{p}^*(M) \\
 & \cong q_! q^*\tilde{p}_!\tilde{p}^*(M)\\
 & \cong q_!q^*p^*p_!(M) \\
 & \cong q_!q^*(M) \cong 0.
\end{align*}
Here we used the definition of $j$ and the commutativity of the left square in the first line, proper base change for the left square in the second line, proper base change for the right square in the third line and  fully faithfulness of $p_!$ and that $M\in \mathrm{ker}(q^*)$ in the last line. \\
4)
For any open immersion $j: U \to X$ we have (for example by the description in terms of idempotent algebras in \Cref{alternative defclosed open inSym}) that the diagonal $\Delta_j : U \to U \times_X U$ is an isomorphism, thus $j$ is a monomorphism in $\Aff$ and hence $(-1)$-truncated.
\end{proof}
\begin{definition}\Cite[Definition 3.3.2]{heyer20246} 
Let $(C,E)$ be a geometric setup such that $C$ admits finite products. A suitable decomposition of $E$ is a pair $I,P \subset E$ such that the conditions 1)-4) of the previous lemma are satisfied.
\end{definition}
Thus the pair $I,P \subset E$ of open immersions and proper maps define a suitable decomposition of $E$ in $\Aff$. 

\begin{lemma}\label{proper maps arre Künneth}
    Any proper map $p: \SpecAn(B) \to \SpecAn(A)$ is Künneth. 
\end{lemma}
\begin{proof}
 Let $g: \SpecAn(C)\to \SpecAn(A)$ be any morphism in $\Aff$. By \Cref{alternative char. of open and proper in aff}, we have an equivalence $D(B)\cong \Mod_{B^{\tri}}D(A)$. Since any base change of a proper map is proper by \Cref{alternative char. of open and proper in aff}, we obtain that $D(C\otimes_AB)\cong \Mod_{g^*(B^{\tri})}D(C)$. On the other hand we have  equivalences 
 \[
 D(C)\otimes_{D(A)}D(B)\cong D(C)\otimes_{D(A)}\Mod_{B^{\tri}}D(A)\cong \Mod_{g^*(B^{\tri})}D(C).
 \]
\end{proof}
\begin{lemma}\label{I and P interact condition c}
Consider a cartesian diagram in $\Aff$
\[\begin{tikzcd}
	{\SpecAn(D)} & {\SpecAn(C)} \\
	\SpecAn(A) & \SpecAn(B)
	\arrow["{j'}", from=1-1, to=1-2]
	\arrow["{p'}", from=1-1, to=2-1]
	\arrow["\ulcorner"{anchor=center, pos=0.125}, draw=none, from=1-1, to=2-2]
	\arrow["p", from=1-2, to=2-2]
	\arrow["j", from=2-1, to=2-2]
\end{tikzcd}\]
with $j\in I$ and $p\in P$. Then the natural map  $j_!p'_*\to p_*j'_! $ is an isomorphism of functors.
\end{lemma}
\begin{proof}
Since $j'$ is an open immersion the functor $j'^*$ is essentially surjective. It thus suffices to show the isomorphism for $M\coloneqq j'^*(N)$, $N\in D(X')$. Let $I \in D(B)$ be the idempotent algebra corresponding to $j$, then since $p$ (and hence $p'$) is proper we have that $C\cong C^{\tri}_{/B}$ the idempotent algebra corresponding to $j'$ is given by by $p^*(I)\cong I \otimes_B C^{\tri}$. We claim that there are the following equivalences:
\begin{align*}
  p_*j'_!j'^*(N)& \cong p_*(\text{fib}(N\to N\otimes_{C}p^*(I)) \\
   &\cong \text{fib}(p_*(N)\to p_*(N)\otimes_{B}I )
\end{align*}
Here we used the definition of $j'_!$ in the first line and that $p_*$ commutes with limits and satisfies proper base change in the second line. On the other hand, we obtain $j_!p'_*j'^*(N)\cong j_!j^*p_*(N)$ by proper base change. Using the definition of $j_!j^*$ we thus obtain 
\begin{align*}
    j_!p'_*j'^*(N)& \cong j_!j^*p_*(N) \\
    & \cong \text{fib}(p_*(N) \to p_*(N)\otimes_B I)
\end{align*}
which shows that the natural map $j_!p'^*\to p_*j_!'$ is an equivalence.
\end{proof}

\begin{proposition}\Cite[lecture 17]{AnSt}\label{affinekünneth}
The pair $(\Aff,E)$ defines a geometric setup and the functor $D: \Aff^{op} \to \PrL$ extends to a $6$-functor formalism  $D: Corr(\Aff,E) \to \PrL$ which is Künneth (and Tannakian if we work with $\text{AnCR}$).
\end{proposition}
\begin{proof}
 To prove that the functor $D$ extends to a $6$-functor formalism we need to show that the assumptions a), b) and c) in \Cite[Proposition 3.3.3]{heyer20246} are verified. The conditions a) and b) are satisfied by definition of the classes $I$ and $P$ and by the base change results for categorical open immersions and proper maps \Cref{alternative char. of open and proper in aff}, \Cref{pbc and stability under pullback of open imersion in aff}. The last condition c) is satisfied by \Cref{I and P interact condition c}. By \Cite[Proposition 3.3.3]{heyer20246} we thus obtain a $6$-functor formalism 
 \[
 D: Corr(\Aff,E) \to \PrL
 \]
 such that for all $p\in P$ the functor $p_!$ is right adjoint to $p^*$ and for all $j\in I$ the functor $j_!$ is left adjoint to $j^*$. \\
 We will now prove that any $!$-able map $f : \SpecAn(A)\to \SpecAn(R)$ is Künneth. By definition $f\cong p \circ  j$ for $p$ a proper map and $j$ an open immersion. By \Cref{proper maps arre Künneth} and \Cref{open immersion and closed immersions are Künneth}, open immersion and proper maps are Künneth, so we conclude by \Cref{dualisbe stable under compposition}. The claim that $D$ is Tannakian is \Cite[Lemma 2.1.3]{camargo2024analytic}.
\end{proof}

We work with the following definition of analytic stacks which is slightly simplified compared to the correct definition in \Cite[Lecture 19]{AnSt}.

\begin{definition}\Cite[Lecture 19]{AnSt}
 An analytic stack is a sheaf $X \in \Shv(\Aff, \Ani)$ for the $!$-topology on $\Aff$.
 
\end{definition}
\begin{remark}
Since the category $\AnR$ is not small one needs to be careful when considering sheaves on it. One way to proceed is to pick a large enough infinite regular cardinal $\kappa$ and to consider $\kappa$-small objects in $\AnR$. 
\end{remark}
By \Cref{affinekünneth} we have a sheafy $6$-functor formalism $D : Corr(\Aff, E) \to \PrL$ on the geometric setup $(\Aff, E)$, with $E$ the class of $!$-able maps. Using the extension result \Cite[Theorem 3.4.11]{heyer20246} we obtain a sheafy $6$-functor formalism
\[
\mathit{D}_{qc} : Corr(\AnSt, \tilde{E}) \to \CAlg(\PrL)
\]
on the category of analytic stacks with $\tilde{E}$ the class of morphisms as in \Cref{E' properties}.

\begin{corollary}
 Let $f: X \to S$ be a morphism of analytic stacks in $\tilde{E}$ and assume that there is a $!$-cover $g: S'\to S$ with $S'\in C$ and $g\in E_0$. Then $f$ is Künneth (and Tannakian if we work with analytic $\mathbb{E}_\infty$-rings).
\end{corollary}
\begin{proof}
    This  follows from \Cref{affinekünneth} together with \Cref{Künnethextension-main corollary}.
\end{proof}

\subsection{Examples of analytic stacks and $!$-covers}
We will give some examples of analytic stacks and $!$-covers.
\begin{example}\Cite[Lecture 19]{AnSt}
For any $A\in \AnR$, the functor 
\[
\SpecAn(A) : F: \AnR \to \Ani , \ B \mapsto \text{Hom}(A,B)
\]
is an analytic stack.
\end{example}
\begin{example}\Cite[Lecture 19]{AnSt}
 Let $S\in \text{Prof}^{\text{light}}$ be a light profinite set, then we associate the analytic stack $\SpecAn(C(S,\mathbb{Z}))$, where $C(S,\mathbb{Z})$ is the ring of $\mathbb{Z}$-valued continuous functions on $S$, with its discrete analytic ring structure. This gives a functor 
 \[
 (-)_{\text{Betti}}: \text{Prof}^{\text{light}} \to \AnSt.
 \]
 One can check that this sends surjections of light profinite sets to $!$-covers and thus gives a functor from the category of compact Hausdorff spaces $\text{LCHaus}$ to analytic stacks
 \[
 (-)_{\text{Betti}}: \text{CHaus} \to \AnSt.
 \]
\end{example}

\begin{example}
    Let us sketch how to realise algebraic stacks in the category $\AnSt$. Let $\text{Ring}$ be the category of animated commutative rings. For $R \in \text{Ring} $ we can associate to it the analytic ring $R^{\tri}\coloneqq (\underline{R},\Mod_{\underline{R}}D(\mathbb{Z}))$), giving a functor 
    \[
    (-)^{\tri}: \text{Ring} \to \AnR .   \]
    This functor sends fpqc covers to $!$-covers, we thus obtain a functor
     \[
    (-)^{\tri}: \Shv_{fpqc}(\text{Ring},\Ani) \to \AnSt.  \]
\end{example}

\begin{example}
 Let $(A,A^+)$ be a sheafy analytic Huber pair, then by \Cite[Proposition 3.34]{andreychev2021pseudocoherent} we can associate to it an analytic ring $(A,A^+)_{\square}$. Denoting by $\text{AnH}$ the category of complete sheafy analytic Huber pairs we thus obtain a functor 
 \[
 (-)_{\square}: \text{AnH} \to \AnSt,  \ (A,A^+) \mapsto \SpecAn((A,A^+)_{\square}).
 \]
 This functor sends open covers to $!$-covers, inducting a functor from the category $\text{AdicSp}$ of analytic adic spaces to analytic stacks
\[
 (-)_{\square}: \text{AdicSp} \to \AnSt.
\]
\end{example}

\subsection{A $p$-adic version of Drinfeld's lemma}\label{Examples of Künneth morphism for analytic stacks}\label{Drinfelds section}
In this subsection we want to give an application of our results which we hope to be of interest for a $p$-adic version of the work of Fargues and Scholze (cf. \Cite{fargues2021geometrization}). A central object in their geometrisation of the $\ell$-adic local Langlands correspondence is the following analytic adic space
\[
Y_{(0,\infty),S}\coloneqq \Spa(W(R^+))\backslash V([\pi]p),
\]
where $W(-)$ denotes the $p$-typical Witt vectors, $S=\Spa(R,R^+)$ an affinoid perfectoid space in characteristic $p$ and $\pi\in R^+$ a pseudo-uniformiser. Let us fix an algebraic closure $k\coloneqq \bar{\mathbb{F}}_p$. One main insight of the work of Anschütz, Le Bras and Mann \cite{anschütz20246functorformalismsolidquasicoherent}, is that the $p$-adic pro-étale cohomology of $S$ can be understood by quasi-coherent cohomology of the quotient of the analytic adic space $Y_{(0,\infty),S}$ by the Frobenius $\varphi$ induced by $R^+$
\[
Y_{(0,\infty),S}/\varphi^{\mathbb{Z}}.
\]
The Drinfeld lemma in the $\ell$-adic setting of Fargues Scholze deals with the $v$-stack \newline $\Div_k^1\coloneqq (\Spa(\mathbb{Q}_p^{\cyc})/\mathbb{Z}_p^\times)/\varphi^{\mathbb{Z}}$ and the full subcategory $D_{dlb}(Y,\Lambda)\subset D_{et}(Y,\Lambda)$ of dualisable objects in the category of étale sheaves with coefficients in $\Lambda\coloneqq \Bar{\Q}_{\ell}$ on a small $v$-stack $Y$. It can be phrased as an equivalence of idempotent complete stable categories 
\[
D_{dlb}(\Div_k^1,\Lambda)\otimes_{D_{dlb}(\Spd(k),\Lambda)}D_{dlb}(W,\Lambda) \cong D_{dlb}(\Div_k^{1}\times_{\Spd(k)}W,\Lambda) 
\]
for $W$ any small $v$-stack over $\Spd(k)$ (cf. \Cite[Chapter IV.7.3]{fargues2021geometrization}). 
Guided by the work of Anschütz, Le Bras and Mann and motivated to obtain a $p$-adic version of Drinfeld's lemma we consider the following analytic stacks

\begin{align*}
    \Div_{\mathbb{Q}_p}^{1,\ct}\coloneqq & ((Y_{(0,\infty), \Q_p^{\cyc,b}})_{\square}/(\Z_p^{\times})_{\text{Betti}})/\varphi^{\Z} \\
\Div_{\mathbb{Q}_p}^{1,\la}\coloneqq & (\lim_{T\mapsto (T+1)^p-1}\mathbb{D}_{\square} \backslash \{0\} )/\Q_p^{\times, \la}.
\end{align*}

By  \Cref{Künnethextension-main theorem}, any morphism $X\to S$ of analytic stacks in $\tilde{E}$ with $S\in \AnSt$ such that $\Delta_S$ is representable and admitting a $!$-cover $U\to S$ with $U \in C$ give examples of Künneth morphisms for the $6$-functor formalism $\mathit{D}_{qc} : Corr(\AnSt, \tilde{E}) \to \PrL$. In particular we obtain the following example.
\begin{corollary}\label{Drinfeld cont und la}
The structural morphism $\Div_{\mathbb{Q}_p}^{1,?} \to \SpecAn(\Q_{p,\square})$  for $?\in \{ \ct, \la \}$ lies in $\tilde{E}$ and is thus Künneth.
\end{corollary}

We end with an example for a morphism which is not Künneth. One way to construct such a morphism is to consider stacks $S$ with non-affine diagonal:

\begin{example}
Let $k$ be an algebraically closed field and consider $E \to \Spec(k)$ an elliptic curve. We consider $E^{\tri}\to \SpecAn(k^{\tri})\coloneqq \ast$ with its discrete induced analytic structure. Then the morphism $f : \ast  \to \ast/E^{\tri}$ is not Künneth. Indeed, note that it lies in $\tilde{E}$ as it has a section. One can show that $f$ is a descendable cover and that $f_*(1)\cong 1$. 
We thus obtain an isomorphism 
\[
D_{qc}(\ast/E^{\tri})\cong \coMod_{f^*f_*(1)}D_{qc}(\ast)\cong D_{qc}(\ast).
\]

In particular, we observe that the morphism
\[
D_{qc}(\ast) \otimes_{D_{qc}(\ast/E^{\tri}) }D_{qc}(\ast) \cong D(\ast) \to D_{qc}(\ast \times_{\ast/E^{\tri} }\ast) \cong D_{qc}(E^{\tri})
\]
is not an isomorphism. Note that $f$ has non-representable diagonal $\Delta_f$ as the pullback along itself gives $\ast \times_{\ast/E^{\tri}} \ast \cong E^{\tri}$ which does not lie in $\Aff$.
\end{example}
\printbibliography
\end{document}